\documentclass[11pt,a4paper,reqno]{amsart}

\usepackage{fullpage}

\usepackage{hyperref}

\usepackage{amsmath,amssymb,mathscinet,dsfont}
\usepackage[all]{xypic}

\numberwithin{equation}{section}

\DeclareMathAlphabet{\mathpzc}{OT1}{pzc}{m}{it}

\mathchardef\mhyph="2D

\theoremstyle{plain}
\newtheorem{theorem}{Theorem}[section]
\newtheorem{proposition}[theorem]{Proposition}

\newtheorem{corollary}[theorem]{Corollary}

\theoremstyle{definition}
\newtheorem{definition}[theorem]{Definition}
\newtheorem{problem}[theorem]{Problem}
\newtheorem{example}[theorem]{Example}

\theoremstyle{remark}
\newtheorem{remark}[theorem]{Remark}

\DeclareMathOperator{\Aut}{Aut}
\DeclareMathOperator{\End}{End}
\DeclareMathOperator{\Hom}{Hom}

\DeclareMathOperator{\Nat}{Nat}
\DeclareMathOperator{\Natb}{\Nat_{\mathrm{b}}}
\DeclareMathOperator{\Irr}{Irr}
\DeclareMathOperator{\Rep}{Rep}

\DeclareMathOperator{\Tr}{Tr}

\DeclareMathOperator{\Ch}{Ch}
\DeclareMathOperator{\tr}{tr}

\newcommand{\C}{{\mathbb C}}

\newcommand{\Z}{{\mathbb Z}}

\newcommand{\T}{{\mathbb T}}

\newcommand{\B}{{\mathcal B}}
\newcommand{\U}{\mathcal U}
\newcommand{\M}{{\mathcal M}}

\newcommand{\cat}[1]{\mathcal{#1}}
\newcommand{\cC}{\cat{C}}
\newcommand{\ccD}{\cat{D}}
\newcommand{\cP}{\cat{P}}

\newcommand{\functor}[1]{\mathcal{#1}}
\newcommand{\fE}{\functor{E}}
\newcommand{\fF}{\functor{F}}
\newcommand{\norm}[1]{\left \| #1 \right \|}

\newcommand{\flr}[1]{\left\lfloor #1 \right\rfloor}

\newcommand{\id}{\mathrm{id}}

\newcommand{\SU}{\mathrm{SU}}

\newcommand{\Id}{\mathrm{Id}}
\newcommand{\Dhat}{\hat{\Delta}}
\newcommand{\un}{{\mathds 1}}
\newcommand{\trhd}{\mathbin{\tilde{\rhd}}}
\newcommand{\alg}{\textrm{alg}}
\newcommand{\Hilb}{\mathrm{Hilb}}
\newcommand{\Hilbf}{\Hilb_{\mathrm{f}}}

\newcommand{\Tens}{\mathpzc{Tens}}
\newcommand{\YD}{\mathpzc{YD}_{\mathrm{brc}}}

\DeclareFontFamily{U}{wasy}{}
\DeclareFontShape{U}{wasy}{m}{n}{ <5> <6> <7> <8> <9> gen * wasy
      <8.5>wasy9
      <10> <10.95> <12> <14.4> <17.28> <20.74> <24.88>wasy10  }{}
\DeclareFontShape{U}{wasy}{b}{n}{ <-10> sub * wasy/m/n
 <10> <10.95> <12> <14.4> <17.28> <20.74> <24.88>wasyb10 }{}
\DeclareFontShape{U}{wasy}{bx}{n}{ <-> sub * wasy/b/n}{}
\DeclareSymbolFont{wasy}{U}{wasy}{m}{n}
\SetSymbolFont{wasy}{bold}{U}{wasy}{b}{n}
\DeclareMathSymbol\ocircle {\mathbin}{wasy}{"23}

\newcommand{\circt}%
{\mathbin{%
\mathchoice
{\ooalign{$\ocircle$\cr\hidewidth\raise-.13ex\hbox{$\scriptstyle\top\mkern2.57mu$}\cr}}
{\ooalign{$\ocircle$\cr\hidewidth\raise-.13ex\hbox{$\scriptstyle\top\mkern2.57mu$}\cr}}
{\ooalign{$\scriptstyle\ocircle$\cr\hidewidth\raise-.18ex\hbox{$\scriptscriptstyle\top\mkern2.3mu$}\cr}}
{\ooalign{$\scriptstyle\ocircle$\cr\hidewidth\raise-.18ex\hbox{$\scriptscriptstyle\top\mkern2.3mu$}\cr}}
}}

\begin{document}

\title{Towards a classification of compact quantum groups of Lie type}

\author[S. Neshveyev]{Sergey Neshveyev}
\email{sergeyn@math.uio.no}
\address{Department of Mathematics, University of Oslo, P.O. Box 1053
Blindern, NO-0316 Oslo, Norway}
\thanks{S.N. was supported by the European Research Council under the European Union's Seventh
Framework Programme (FP/2007-2013)/ERC Grant Agreement no. 307663
}

\author[M. Yamashita]{Makoto Yamashita}
\email{yamashita.makoto@ocha.ac.jp}
\address{Department of Mathematics, Ochanomizu University, Otsuka
2-1-1, Bunkyo, 112-8610 Tokyo, Japan}
\thanks{M.Y. was supported by JSPS KAKENHI Grant Number 25800058, and partially by the Danish National Research Foundation through the Centre for Symmetry and Deformation (DNRF92)}

\thanks{This survey is to appear in the volume ``The Abel Symposium 2015 -- Operator Algebras and Applications'', the final publication will be available at \url{http://link.springer.com}. 
}

\date{March 16, 2015}

\begin{abstract}
This is a survey of recent results on classification of compact quantum groups of Lie type, by which we mean quantum groups with the same fusion rules and dimensions of representations as for a compact connected Lie group~$G$. The classification is based on a categorical duality for quantum group actions recently developed by De Commer and the authors in the spirit of Woronowicz's Tannaka--Krein duality theorem. The duality establishes a correspondence between the actions of a compact quantum group $H$ on unital C$^*$-algebras and the module categories over its representation category $\Rep H$. This is further refined to a correspondence between the braided-commutative Yetter--Drinfeld $H$-algebras and the tensor functors from~$\Rep H$. Combined with the more analytical theory of Poisson boundaries, this leads to a classification of dimension-preserving fiber functors on the representation category of any coamenable compact quantum group in terms of its maximal Kac quantum subgroup, which is the maximal torus for the $q$-deformation of $G$ if~$q\ne1$. Together with earlier results on autoequivalences of the categories $\Rep G_q$, this allows us to classify up to isomorphism a large class of quantum groups of $G$-type for compact connected simple Lie groups $G$. In the case of $G=\SU(n)$ this class exhausts all non-Kac quantum groups.
\end{abstract}

\maketitle

\section{Introduction}
\label{sec:intro}

The theory of quantum groups has two origins.  One is the algebraic approach motivated by the quantum inverse scattering method and initiated by the discovery of quantized universal enveloping algebras by Drinfeld~\cite{MR802128} and Jimbo~\cite{MR797001}. The other is the operator algebraic approach developed by Woronowicz~\cite{MR901157}, which stands on the philosophy of treating noncommutative C$^*$-algebras as functions on `noncommutative spaces', or `pseudospaces' in Woronowicz's terminology. In both frameworks, deformations of $\SU(2)$ and, more generally, of the compact simple Lie groups appear as the most fundamental and motivating examples.

A natural question is to classify all such deformations within a reasonable system of axioms.  At the infinitesimal level, this was already settled in a series of papers by Drinfeld and his collaborators during the 1980s. If one assumes compatibility with the compact form, the classification takes a particularly elegant form~\cite{MR714225,MR1047964,MR1017083}: if $\mathfrak{g}$ is a complex simple Lie algebra and~$G$ is the simply connected compact Lie group corresponding to $\mathfrak{g}$, any infinitesimal deformation of the universal enveloping algebra $\U(\mathfrak{g})$ as a quasitriangular $*$-Hopf algebra corresponds to a Poisson--Lie group structure on $G$, and those are parametrized, up to isomorphisms and rescalings, by invariant Poisson structures on a maximal torus.  Moreover, these deformations make perfect sense in Woronowicz's framework and give rise to strict deformation quantization of the corresponding Poisson structures in the sense of Rieffel~\cite{MR1292010}.

Extending this classification from infinitesimal to analytical level, or in other words, going from formal to strict deformation quantization, is a nontrivial task.  Since there is no analogue of Kontsevich's classification theorem in the analytic setting, we cannot expect such a clean result as in the formal case without imposing first an additional symmetry. One possible idea is to require the deformations to have the same combinatorial structure of representation theory, meaning the fusion rules and dimensions of representations, as in the classical case. This implies that we can take the coalgebra of matrix coefficients independent of the deformation parameter, and then try to find a suitable algebra structure on it. Moreover, if we want to preserve quasitriangularity, the Faddeev--Reshetikhin--Takhtadzhyan method~\cite{MR1015339} reduces the problem to finding $R$-matrices with certain algebraic symmetry. Still, carrying out such a classification directly is not an easy task, and this has been only worked out for the quantum groups with representation theory of $\SU(2)$~\cite{MR890482} and of $\SU(3)$~\cite{MR1673475,MR2106933}.

Our approach to this problem consists of a combination of cohomological and analytical methods for tensor categories and related constructions for quantum groups. Suppose we want to classify all quantum groups which have the same fusion rules and dimensions of representations as a compact group $G$. In view of Woronowicz's Tannaka--Krein duality theorem, the problem can be divided into three parts:
\begin{enumerate}
\item[--] classification of rigid C$^*$-tensor categories $\cC$ with fusion rules of $G$;
\item[--] classification of monoidal autoequivalences of $\cC$;
\item[--] classification of unitary fiber functors $\cC\to\Hilbf$ inducing the classical dimension function on the representation ring of $G$.
\end{enumerate}
We will mainly concentrate on the third problem, where the most recent advances are. Let us explain the strategy in more detail, in roughly the same order as we proceed from Section~\ref{sec:cat-duality-action} to Section~\ref{sec:q-grp-Lie-type}.

The first step is to establish a duality between the category of unital $H$-C$^*$-algebras and that of $(\Rep H)$-module categories for compact quantum groups $H$ \cite{MR3121622,DOI:10.17879/58269764511}. Given an $H$-C$^*$-algebra $B$, which represents a noncommutative $H$-space $X$, we take the category of $H$-equivariant Hermitian vector bundles on $X$, that is, of finitely generated projective $H$-equivariant Hilbert $B$-modules. We can consider the tensor product of such modules and finite dimensional unitary representations of~$H$, which leads to the structure of a $(\Rep H)$-module category. An analogue of the reconstruction procedure in Woronowicz's Tannaka--Krein theorem implies that $B$ can be recovered from this module category and the distinguished object in it represented by $B$ itself.

Further pursuing this duality, we have the following correspondence for subclasses of $H$-C$^*$-algebras and $(\Rep H)$-module categories~\cite{MR3291643}. Among the $H$-C$^*$-algebras we consider the braided-commutative Yetter--Drinfeld algebras, while among the module categories we take the C$^*$-tensor categories with module structure induced by a tensor functor from $\Rep H$. A motivating example for this duality is the coideal of the function algebra coming from a quantum subgroup $K$ of $H$, on one side, and the forgetful functor $\Rep H \to \Rep K$, on the other. In this formulation the quantum subgroup coideals precisely correspond to the factorizations $\Rep H \to \cC \to \Hilbf$ of the canonical fiber functor on $\Rep H$ through some C$^*$-tensor category~$\cC$.

One of our discoveries is that the noncommutative Poisson boundary of the discrete dual of~$H$, and its counterpart on the categorical side, fits very well into this scheme~\cite{arXiv:1405.6572}. The noncommutative Poisson boundary, modeled on the classical theory for discrete groups initiated by Furstenberg~\cite{MR0352328},  was introduced by Izumi~\cite{MR1916370} to understand the lack of minimality for infinite tensor product actions of non-Kac quantum groups. This theory requires the operator algebraic framework in an essential way, and that is why we need to work with compact quantum groups instead of, for example, cosemisimple Hopf algebras. We find that the categorical Poisson boundary has a universality property for what we call \emph{amenable tensor functors}. A concrete implication is that if $G_q$ is the Drinfeld--Jimbo $q$-deformation of a compact connected semisimple Lie group $G$ with a maximal torus $T$, then for $q\ne1$ the forgetful functor $\Rep G_q\to\Rep T$ is a universal tensor functor defining the classical dimension function on $\Rep G$. This shows that the undeformed classical torus $T < G_q$ can be detected already at the categorical level. This also implies that the dimension-preserving fiber functors on $\Rep G$ are parametrized by the invariant Poisson structures on the torus, or the Poisson--Lie group structures of $G$, as expected from the case of infinitesimal deformations.

Returning to the first problem of classifying rigid C$^*$-tensor categories with fusion rules of~$G$, there is a simple class of examples of such categories in addition to $\Rep G_q$.
We have a natural grading of $\Rep G_q$ by the Pontryagin dual of the center of~$G$, or more intrinsically, by the chain group of $\Rep G_q$. Then any $\T$-valued $3$-cocycle on the chain group defines an associator, and we obtain a twisted tensor category, which has the same fusion rules as $G$. Note that this is analogous to, but much easier than, the famous Knizhnik--Zamolodchikov associator constructed by Drinfeld, which relates $\Rep G$ and $\Rep G_q$. What we wrote in the previous paragraph about the categories $\Rep G_q$ applies equally well to these twisted categories and allows us to classify dimension-preserving fiber functors on them. In the case of~$\SU(n)$, up to monoidal equivalence, these exhaust the tensor categories with fusion rules of $\SU(n)$~\cite{MR1237835,MR3266525}. It seems reasonable to expect that a similar result it true, or at least close to be true, for other simple Lie groups. Finally, the second problem of classifying autoequivalences of such categories has been essentially solved in~\cite{MR2844801,MR2959039}.

To summarize our results, for any compact connected simple Lie group $G$, we classify up to isomorphism all compact quantum groups with the same fusion rules and dimensions of representations as for $G$ which, moreover, have representation categories $\Rep G_q$ for $q\ne1$ or twists of such categories by $3$-cocycles on the dual of the center of $G$, see Theorems~\ref{thm:3coc-2pseudo} and~\ref{thm:simple-G-isom-class}. Whether such categories exhaust all categories with fusion rules of $G$ beyond the case of $\SU(n)$, as well as whether one can say something precise for $q=1$, are the main remaining open questions.

\section{Monoidal categories}

Throughout the exposition we fix a universe and assume that all categories are \emph{small}~\cite{MR1712872}. See~\cite{MR3204665} for the details on the following notions.

\subsection{C\texorpdfstring{$^*$}{*}-categories}

A \emph{C$^*$-category} is a category $\cC$ with morphism sets
$\cC(U, V)$ that are complex Banach spaces, endowed with a complex conjugate
involution $\cC(U, V) \to \cC(V, U)$, $T \mapsto T^*$, satisfying the
C$^*$-identity
$$\norm{T^* T} = \norm{T}^2 = \norm{T T^*},$$
and such that the composition of morphisms is bilinear and $\|ST\|\le\|S\|\,\|T\|$.
Unless said otherwise, we always assume that $\cC$ is
closed under finite direct sums and subobjects. The latter means that any
idempotent in the endomorphism ring $\cC(X) = \cC(X, X)$ comes
from a direct summand of~$X$. The existence of finite direct sums guarantees that, for any $T\colon X\to Y$, the morphism $T^*T$ is positive as an element in the C$^*$-algebra $\cC(X)$, which otherwise we would have to add as an additional axiom.

An object is called \emph{simple} if its endomorphism ring is isomorphic to $\C$, and a C$^*$-category is said to be \emph{semisimple} if any object is isomorphic to a finite direct sum of simple ones. We then denote the isomorphism classes of simple
objects by~$\Irr \cC$ and assume that this is an at most countable
set. Many results admit formulations which do not require this
assumption, and can be proved by considering subcategories generated by
countable sets of simple objects, but we leave this matter to the
interested reader.

A \emph{unitary functor}, or a \emph{C$^*$-functor}, is a linear
functor of C$^*$-categories $\fF\colon \cC \to \cC'$ satisfying $\fF(T^*)
= \fF(T)^*$.

A few times we will need to perform the following operation: starting
from a C$^*$-category~$\cC$, we replace the morphisms sets by some
larger multiplicative system $\ccD(X, Y)$ naturally containing the original $\cC(X,
Y)$.  Then we perform the \emph{idempotent completion} to construct a
new category~$\ccD$. That is, we regard the projections $p
\in \ccD(X)$ as objects in the new category, and take $q \ccD(X,
Y) p$ as the morphism set from the object represented by $p
\in \ccD(X)$ to the one by $q \in \ccD(Y)$.  Then the
embeddings $\cC(X, Y) \to \ccD(X, Y)$ can be considered as a
C$^*$-functor $\cC \to \ccD$.

\subsection{C\texorpdfstring{$^*$}{*}-tensor categories}\label{sec:c-star-cat}

A \emph{C$^*$-tensor category} is a C$^*$-category endowed with a
unitary bifunctor $$\mathord{\otimes} \colon \cC \times \cC \to \cC,$$ a
distinguished object $\un \in \cC$, and natural unitary isomorphisms
\begin{align*}
\un \otimes U &\cong U \cong U \otimes \un,& \Phi({U,
V, W})&\colon (U \otimes V) \otimes W \to U \otimes (V \otimes W)
\end{align*}
satisfying certain compatibility conditions. If these isomorphisms can be taken to be the identity morphisms, then $\cC$ is said to be \emph{strict}.

We denote by $\Hilbf$ a category of finite dimensional Hilbert spaces with a strict model of tensor product (see, for example,~\cite{MR1870871} or~\cite[p.~37]{MR3204665} for concrete realizations). 

A \emph{unitary tensor functor}, or a \emph{C$^*$-tensor functor},
between two C$^*$-tensor categories $\cC$ and $\cC'$ is given by a
triple $(\fF_0, \fF, \fF_2)$, where $\fF$ is a C$^*$-functor $\cC \to \cC'$,
$\fF_0$ is a unitary isomorphism $\un_{\cC'} \to \fF(\un_\cC)$, and $\fF_2$
is a collection of natural unitary isomorphisms $\fF(U) \otimes \fF(V) \to \fF(U \otimes
V)$ compatible with the structure morphisms of $\cC$ and
$\cC'$.

As a rule, we denote tensor functors by just one symbol $\fF$. A tensor functor between strict C$^*$-tensor categories is said to be \emph{strict} if $\fF_0$ and $\fF_2$ are the identity morphisms. The composition of tensor functors $\fF$ and $\fF'$ is defined by taking the usual composition of C$^*$-functors and setting $(\fF \fF')_2 = \fF(\fF'_2) \fF_2$ and $\fF\fF'=\fF(\fF_0')\fF_0$. A \emph{natural transformation} of tensor functors is a natural transformation in the usual sense which is also compatible with the isomorphisms $\fF_2$ and $\fF_0$.

When the C$^*$-functor part of $\fF$ is an equivalence of categories, $\fF$ is said to be a \emph{unitary monoidal equivalence}. Analogously to the case of ordinary functors, unitary monoidal equivalences can be inverted (up to a natural unitary monoidal isomorphism) as unitary tensor functors, 
so for any $\fF$ as above there is another unitary monoidal equivalence $\fF' \colon \cC' \to \cC$ such that $\fF \fF'$ and $\fF' \fF$ are naturally unitarily monoidally isomorphic to the identity functors of the respective categories. If the target and source categories are the same, such a tensor functor~$\fF$ is called an \emph{autoequivalence}, and we denote by $\Aut^\otimes(\cC)$ the group of autoequivalences of~$\cC$ considered up to natural unitary monoidal isomorphisms. Let us also note that a version of Mac Lane's coherence theorem~\cite[Chapter~XI]{MR1712872} says that any C$^*$-tensor category is unitarily monoidally equivalent to a strict one.

\smallskip

When $\cC$ is a C$^*$-tensor category and $U \in \cC$, an
object $V$ is said to be a \emph{dual object} of $U$ if there are
morphisms $R \in \cC(\un, V \otimes U)$ and $\bar{R} \in \cC(\un, U
\otimes V)$ satisfying the conjugate equations
\begin{align*} (\iota_{V} \otimes \bar{R}^*) \Phi (R \otimes \iota_{V}) &=
\iota_{V},& (\iota_{U} \otimes R^*) \Phi (\bar{R} \otimes \iota_{U}) &=
\iota_{U},
\end{align*}
where we assumed for simplicity that the unit is strict, so that $\un\otimes U=U\otimes\un=U$.
If any object in $\cC$ admits a dual, $\cC$ is said to be
\emph{rigid} and we denote a choice of a dual of $U \in \cC$
by~$\bar{U}$.

Any rigid C$^*$-tensor category with simple unit has
finite dimensional morphism spaces and hence is automatically
semisimple. The quantity
$$
d^\cC(U) =\min_{(R, \bar{R})} \norm{R} \norm{\bar{R}},
$$
where $(R, \bar{R})$
runs through the set of solutions of conjugate equations as above, is called the \emph{intrinsic dimension} of $U$. We
omit the superscript $\cC$ when there is no danger of confusion.  A
solution $(R, \bar{R})$ of the conjugate equations for $U$ is called
\emph{standard}~if
$$
\|R\|=\|\bar R\|=d(U)^{1/2}.
$$
Solutions of the conjugate equations for $U$ are unique up to the
transformations $$(R, \bar{R}) \mapsto ((T^* \otimes \iota) R, (\iota
\otimes T^{-1}) \bar{R}).$$ Furthermore, if $(R,\bar R)$ is standard,
then such a transformation defines a standard solution if and only if
$T$ is unitary.

In a rigid C$^*$-tensor category $\cC$ we often fix standard solutions
$(R_U,\bar R_U)$ of the conjugate equations for every object $U$. Then
$\cC$ becomes \emph{spherical} in the sense that one has the equality
$R_U^* (\iota \otimes T) R_U = \bar{R}_U^* (T \otimes \iota)
\bar{R}_U$ for any $T \in \cC(U)$.  The normalized linear
functional
$$
\tr_U(T) = \frac{1}{d(U)} R_U^* (\iota \otimes T)
R_U=\frac{1}{d(U)}\bar{R}_U^* (T \otimes \iota) \bar{R}_U
$$
is a tracial state on the finite dimensional C$^*$-algebra
$\cC(U)$, called the \emph{normalized categorical trace}. It is independent of the choice of a standard
solution.


\smallskip

For any semisimple C$^*$-tensor category we denote by $K^+(\cC)$ the Grothendieck semiring of~$\cC$. As an additive semigroup it is generated by the isomorphism classes~$[U]$ of objects in $\cC$ and satisfies the relations $[U\oplus V]=[U]+[V]$, so it is the free commutative semigroup with generators $[U]\in\Irr \cC$. The product is defined by $[U]\,[V]=[U\otimes V]$.

Let us fix representatives $\{U_s\}_{s\in\Irr \cC}$ of the isomorphism classes of simple objects of $\cC$. We will often use the subindex $s$ to denote various construction related to~$U_s$. For every object $X\in\cC$, denote by $\Gamma_X=(a^X_{st})_{s,t\in\Irr \cC}$ the matrix describing the multiplication by $[X]$,
$$
[X]\,[U_t]=\sum_s a^X_{st}[U_s],
$$
so $a^X_{st}=\dim\cC(U_s,X\otimes U_t)$. It is not difficult to see that $\|\Gamma_X\|\le d^\cC(X)$. Moreover, the same is true for any \emph{dimension function} $d$ in place of $d^\cC$, by which one means a unital semiring homomorphism $d\colon K^+(\cC)\to[0,+\infty)$ such that $d([X])=d([\bar X])$ for all~$X$.

\smallskip

Another algebraic structure naturally associated with a semisimple C$^*$-tensor category $\cC$ is the \emph{chain group} $\Ch(\cC)$. It is the group with generators $g_s$, $s \in \Irr \cC$, satisfying the relations $g_r = g_sg_t$ whenever $U_r$ embeds into $U_s \otimes U_t$. It is closely related to the notion of \emph{grading} on a category. For a discrete group $\Gamma$, we say that $\cC$ is \emph{graded over} $\Gamma$ if we are given full subcategories $(\cC_g)_{g\in\Gamma}$ such that $\un\in\cC_e$, any object of~$\cC$ is isomorphic in an essentially unique way to a finite direct sum~$\oplus_g X_g$, with $X_g \in \cC_g$, and $X \otimes Y$ is isomorphic to an object in $\cC_{g h}$ if $X \in \cC_g$ and $Y \in \cC_h$. The chain group $\Ch(\cC)$ defines a grading on $\cC$: for every $g\in\Ch(\cC)$, the subcategory~$\cC_g$ consists of direct sums of simple objects $U_s$ such that $g_s=g$. Any other grading over a group~$\Gamma$ defines a canonical homomorphism $\Ch(\cC)\to\Gamma$, so the chain group is a universal group over which $\cC$ is graded.

\section{Compact quantum groups}

Compact quantum groups, or compact matrix pseudogroups as introduced by Woro\-now\-icz~\cite{MR901157}, are our main object of interest. We again follow the presentation of~\cite{MR3204665}.

\subsection{Tannaka--Krein duality}

A \emph{compact quantum} group $G$ is represented by a unital C$^*$-algebra $C(G)$ equipped with a unital $*$-homomorphism $\Delta \colon C(G)\to C(G)\otimes C(G)$ satisfying
\begin{itemize}
\item coassociativity: $(\Delta \otimes \iota) \Delta = (\iota \otimes \Delta) \Delta$,
\item cancellation properties: the linear spans of $$(C(G) \otimes 1) \Delta(C(G))\ \ \text{and}\ \ (1 \otimes C(G)) \Delta(C(G))$$
are dense in $C(G) \otimes C(G)$.
\end{itemize}
There is a unique state $h$ on $C(G)$ satisfying $(h \otimes \iota) \Delta = h$ (and/or $(\iota \otimes h) \Delta = h$) called the \emph{Haar state}.  If $h$ is faithful, $C(G)$ is called the \emph{reduced function algebra of~$G$}, or $G$ is called a \emph{reduced} compact quantum group, and we are mainly interested in such cases. By taking the image of $C(G)$ in the GNS representation of $h$, we can always work with a reduced model.

A \emph{finite dimensional unitary representation} of $G$ is a unitary element $U\in B(H_U)\otimes C(G)$ such that $(\iota\otimes\Delta)(U)=U_{12}U_{13}$, where~$H_U$ is a finite dimensional Hilbert space. The \emph{intertwiners} between two representations $U$ and $V$ are the linear maps $T$ from $H_U$ to $H_V$ satisfying $V (T \otimes 1) = (T \otimes 1) U$. The tensor product of two representations $U$ and $V$ is defined by $U_{13}V_{23}$ and denoted by $U\circt V$. This way, the category $\Rep G$ of finite dimensional unitary representations with intertwiners as morphisms and with tensor product $\circt$ becomes a semisimple C$^*$-tensor category.

%
When $\omega$ is in the (pre)dual $B(H_U)_* \cong \bar{H}_U \otimes H_U$, the element $(\omega \otimes \iota) (U) \in C(G)$ is called the matrix coefficient associated with $\omega$. The dense $*$-subalgebra of $C(G)$ spanned by the matrix coefficients of finite dimensional representations is called the \emph{regular algebra of $G$}, and is denoted by $\C[G]$. This space is closed under the coproduct, and becomes a Hopf $*$-algebra. Its antipode is characterized by
$$
(\iota \otimes S)(U) = U^*\ \ \text{for}\ \ U \in \Rep G.
$$

Let us put $\U(G)=\C[G]^*$. This space has the structure of a $*$-algebra, defined by duality from the Hopf $*$-algebra $(\C[G],\Delta)$. Every finite dimensional unitary representation $U$ of $G$ defines a $*$-representation $\pi_U$ of $\U(G)$ on $H_U$ by $\pi_U(\omega)=(\iota\otimes\omega)(U)$. When convenient, we omit $\pi_U$ and write $\omega\xi$ instead of $\pi_U(\omega)\xi$ for $\xi\in H_U$.

An important ingredient of the duality for unitary representations of $G$ is the Woronowicz character $f_1\in\U(G)$, which we denote by~$\rho$. Namely, the element $\rho\in\U(G)$ is uniquely determined by the properties that it is positive, invertible, $\Tr(\pi_U(\rho))=\Tr(\pi_U(\rho^{-1}))$, and
$$
\bar U=(j(\rho)^{1/2}\otimes1)(j\otimes\iota)(U^*)(j(\rho)^{-1/2}\otimes1)\in B(\bar H_U)\otimes\C[G],
$$
is unitary, where $U$ is any finite dimensional unitary representation and $j$ denotes the canonical $*$-anti-isomorphism $B(H_U)\cong B(\bar H_U)$ defined by $j(T)\bar\xi=\overline{T^*\xi}$. This element $\bar{U}$, which is again a unitary representation, is called the \emph{conjugate representation} of $U$.
The representation $\bar U$ is dual to $U$ in the sense of Section~\ref{sec:c-star-cat}, and a convenient choice of solutions of the conjugate equations for $U$ is given~by
$$
R_U(1)=\sum_i\bar\xi_i\otimes\rho^{-1/2}\xi_i\ \ \text{and}\ \ \bar R_U(1)=\sum_i\rho^{1/2}\xi_i\otimes\bar\xi_i,
$$
where $\{\xi_i\}_i$ is an orthonormal basis in $H_U$. This solution is standard. In particular, the intrinsic dimension function on $\Rep G$ coincides with the quantum dimension
$$
\dim_qU=\Tr(\pi_U(\rho)).
$$

Therefore $\Rep G$ is a rigid C$^*$-tensor category. The forgetful functor $U \mapsto H_U$ defines a strict unitary tensor functor into $\Hilbf$ called the \emph{canonical fiber functor of~$G$}.  Woronowicz's Tannaka--Krein duality theorem recovers the $*$-Hopf algebra~$\C[G]$ from these categorical data.

\begin{theorem}[\cite{MR943923}]\label{thm:Wor-Tan-Kre-thm}
Let $\cC$ be a rigid C$^*$-tensor category with simple unit, and $\fF\colon \cC \to \Hilbf$ be a unitary tensor functor. Then there exist a compact quantum group~$G$, a unitary monoidal equivalence $\fE\colon \Rep G \to \cC$, and a natural unitary monoidal isomorphism from the canonical fiber functor of $G$ to $\fF \fE$. Moreover, the Hopf $*$-algebra $\C[G]$ is determined uniquely up to isomorphism.
\end{theorem}

The key idea is that if $\cC\subset\Hilbf$ and $\fF$ is the embedding functor, then, with representatives~$\{U_s\}_s$ of $\Irr \cC$ and $H_s = \fF(U_s)$, the coalgebra
$$
\bigoplus_{s \in \Irr \cC} B(H_s)_* \cong \bigoplus_{s \in \Irr \cC} \bar{H}_s \otimes H_s
$$
admits an associative product induced by irreducible decompositions of tensor products, which makes it into a bialgebra, analogously to the description of the product of matrix coefficients for usual  compact groups. Moreover,  standard solutions of the conjugate equations determine the involution by the formula
\begin{equation}\label{eq:star-on-CG}
\bar{H}_s \otimes H_s \ni \bar{\xi} \otimes \eta \mapsto \overline{(\iota \otimes \xi^*) R_{s}} \otimes (\eta^* \otimes \iota) \bar{R}_{s} 
\in \bar{H}_{\bar{s}} \otimes H_{\bar{s}}.
\end{equation}
This way one obtains the $*$-bialgebra $\C[G]$.



\smallskip

The representation semiring of $G$ is defined as $R^+(G)=K^+(\Rep G)$. A possible way of saying when a compact quantum group is a deformation of a genuine group is as follows.

\begin{definition}\label{def:type}
Given a compact group $H$, we say that a compact quantum group $G$ is of \emph{$H$-type} if there exists a semiring isomorphism $R^+(G)\cong R^+(H)$ preserving the (classical) dimensions of representations.
\end{definition}

This definition is essentially due to Woronowicz~\cite{MR943923}. In~\cite{MR1734250} the Hopf $*$-algebras $\C[G]$ for compact quantum groups $G$ of $H$-type are called \emph{dimension-preserving $R^+$-deformations} of~$\C[H]$.

\begin{problem}[cf.~\cite{MR943923}]
\label{prob:classification}
Given a compact connected Lie group $H$, classify the compact quantum groups of $H$-type.
\end{problem}

Our aim is to develop a general method to attack this problem, as outlined in Section~\ref{sec:intro}.

\subsection{Cohomology of the discrete dual}
\label{sec:cohom-discr-dual}

Deformation problems for compact quantum groups are controlled by a cohomology theory of the dual discrete quantum groups, which plays a central role in our considerations. We again refer the reader to~\cite{MR3204665} for a more thorough discussion.

\smallskip
The $*$-algebra $\U(G)$ can be identified with  $\prod_{s \in \Irr G} B(H_s)$ (the algebraic direct product) using the correspondence
$$
\C[G]^* \ni \omega \leftrightarrow (\pi_{U_s}(\omega))_s \in \prod_{s \in \Irr G} B(H_s).
$$
More generally, we can consider the algebra
$$
\U(G^k) =(\C[G]^{\otimes k})^*\cong \prod_{s_1,\cdots, s_k \in \Irr G} B(H_{s_1}) \otimes \cdots \otimes B(H_{s_k}),
$$
and interpret it as the space of (possibly unbounded) $k$-point functions on the ``discrete dual'' quantum group $\hat{G}$. This is, for example, how the discrete quantum groups are defined in~\cite{MR1378538}. If $G$ is a genuine commutative compact group, this agrees with the usual notion of functions on the $k$-th power of the Pontryagin dual group $\hat{G}$.

The spaces $\U(G^k)$ can be regarded as the components of the standard complex for the group cohomology with multiplicative scalar coefficients as follows. We call the invertible elements of~$\U(G^k)$ the \emph{$k$-cochains} on $\hat{G}$. Given a $k$-cochain $\omega$, we put
\begin{gather*}
\partial^0(\omega) = 1 \otimes \omega, \quad \partial^{k+1}(\omega) = \omega \otimes 1,\\
\partial^j(\omega) =  \hat{\Delta}_j (\omega) \quad (\hat{\Delta} \text{ applied to the $j$-th position}),
\end{gather*}
(which are all in $\U(G^{k+1})$), and call
$$
\partial(\omega) = \partial^0(\omega) \partial^2(\omega) \cdots \partial^{2\flr{\frac{k+1}2}}(\omega)\partial^1(\omega^{-1}) \partial^3(\omega^{-1}) \cdots \partial^{2\flr{\frac{k}2}+1}(\omega^{-1})
$$
the \emph{coboundary} of $\omega$. When $\partial(\omega) = 1$, $\omega$ is said to be a \emph{$k$-cocycle}. We denote the set of $k$-cocycles by $Z^k(\hat G;\C^\times)$. Two $k$-cocycles $\chi$ and $\chi'$ are said to be \emph{cohomologous}~if
$$
\chi' = \partial^0(\omega) \partial^2(\omega) \cdots \partial^{2\flr{\frac{k}2}}(\omega) \chi \partial^1(\omega^{-1}) \partial^3(\omega^{-1}) \cdots \partial^{2\flr{\frac{k-1}2}+1}(\omega^{-1})
$$
holds for some $(k-1)$-cochain $\omega$.

In general the relation of being cohomologous is not even symmetric. But if it happens to be an equivalence relation, the set of equivalence classes of $k$-cocycles is denoted by $H^k(\hat{G}; \C^\times)$ and called the \emph{$k$-cohomology} of $\hat{G}$. When one requires all the ingredients to be unitary instead of invertible, the corresponding set is denoted by $H^k(\hat{G}; \T)$. If a $k$-cochain $\omega$ commutes with the image of $\hat{\Delta}^{k-1}$ (defined inductively as $(\hat{\Delta}^{k-2} \otimes \iota)\Dhat$), we say that $\omega$ is \emph{invariant}. Invariant cocycles are also called \emph{lazy} in the algebraic literature. If we consider only invariant cocycles (with coboundaries defined also using only invariant cochains), we get sets $H^k_{G}(\hat{G}; \C^\times)$, and $H^k_{G}(\hat{G}; \T)$ for the unitary case, whenever these sets are well-defined.

When $G$ is a commutative group, $H^k_G(\hat{G}; \T) = H^k(\hat{G}; \T)$ agrees with the usual group cohomology of the Pontryagin dual. We also note that, in general, if $H$ is a quantum subgroup of~$G$, the natural inclusion $\U(H^k) \to \U(G^k)$ induces maps from the cohomologies of $\hat{H}$ to those of~$\hat{G}$.

\medskip
We are mainly interested in the cohomology in low degrees ($k \le 3$), which have direct connections with various aspects of $\Rep G$. Let us briefly summarize these connections.

\smallskip

The set $H^1(\hat{G}; \T)=Z^1(\hat{G}; \T)$ consists of unitary group-like elements in~$\U(G)$. If~$G$ is a genuine compact group, then any such group-like element arises from an element of $G$, so we have canonical isomorphisms
$$H^1(\hat{G}; \T)\cong G\ \ \text{and}\ \ H^1_G(\hat{G}; \T)\cong Z(G).$$  For arbitrary $G$, a categorical interpretation of $H^1(\hat{G}; \T)$ is that this is the group of unitary monoidal automorphisms of the canonical fiber functor $\fF\colon \Rep G\to\Hilbf$, while the group $H^1_G(\hat{G}; \T)$  is the group of unitary monoidal automorphisms of the identity functor on $\Rep G$.

\smallskip
Next, a unitary $2$-cochain $x \in \U(G^2)$ is a $2$-cocycle if and only if it satisfies
$$
(x \otimes 1) (\hat{\Delta} \otimes \iota)(x) = (1 \otimes x)(\iota \otimes \hat{\Delta})(x).
$$
If $c$ is a unitary element in the center of $\U(G)$ (an invariant unitary $1$-cochain), then its coboundary is $(c \otimes c) \hat{\Delta}(c^{-1})$. The set of invariant unitary $2$-cocycles forms a group under multiplication, and the coboundaries form a subgroup. Thus, $H^2_G(\hat{G}; \T)$ becomes a group, called the \emph{invariant $2$-cohomology group} of $\hat G$.

If $x$ is an invariant unitary $2$-cocycle, the multiplication by $x^{-1}$ on $H_U \otimes H_V$ can be considered as a unitary endomorphism of $U \circt V$ in $\Rep G$. Such endomorphisms form a natural unitary transformation of the bifunctor $\circt$ into itself. The cocycle condition corresponds to the fact that this transformation is a monoidal autoequivalence of $\Rep G$. Up to natural unitary monoidal isomorphisms, any autoequivalence of $\Rep G$ fixing the irreducible classes can be obtained in this way. Moreover, the cohomology relation of cocycles corresponds to the natural unitary monoidal isomorphism of autoequivalences. Thus $H^2_G(\hat{G}; \T)$ can be considered as the normal subgroup of $\Aut^\otimes(\Rep G)$ consisting of (isomorphism classes of) autoequivalences that preserve the isomorphism classes of objects. Without the unitarity, $H^2_G(\hat{G}; \C^\times)$ corresponds to a subgroup of monoidal autoequivalences of $\Rep G$ as a tensor category over $\C$.

\smallskip

Let $x$ be an arbitrary unitary $2$-cocycle on $\hat{G}$, invariant or not.  Then the triple $\fF_x = (\id_\C, U \mapsto H_U, x^{-1})$ defines a new unitary tensor functor $\Rep G \to \Hilbf$. By Theorem~\ref{thm:Wor-Tan-Kre-thm}, $\fF_x$ can be considered as the canonical fiber functor of another compact quantum group $G_x$ satisfying $\Rep G = \Rep G_x$. Concretely, $\U(G_x)$ coincides with $\U(G)$ as a $*$-algebra, but is endowed with the modified coproduct $\hat{\Delta}_x(T) = x \hat{\Delta}(T) x^{-1}$.  By duality, $\C[G_x]$ is the same coalgebra as $\C[G]$, but has a modified $*$-algebra structure (pre)dual to $(\U(G),\hat{\Delta}_x)$.

Up to natural unitary monoidal isomorphisms, the functors $\fF_x$ exhaust all unitary tensor functors $\fF'\colon\Rep G\to\Hilbf$ satisfying $\dim \fF'(U) = \dim H_U$. Moreover, if~$T$ is a unitary element in $\U(G)$, then $x_T = (T \otimes T) x \hat{\Delta}(T^{-1})$ defines another unitary fiber functor which is naturally unitarily monoidally isomorphic to $\fF_x$. Therefore the set $H^2(\hat{G}; \T)$ gives a parametrization of the set of isomorphism classes of dimension-preserving unitary fiber functors $\Rep G\to\Hilbf$.

The group $H^2_G(\hat{G}; \T)$ acts on the set $H^2(\hat{G}; \T)$ by multiplication on the right. This corresponds to the restriction of the obvious right action of $\Aut^\otimes(\Rep G)$ on the natural unitary monoidal isomorphism classes of unitary fiber functors $\Rep G \to \Hilbf$. Note that by Theorem~\ref{thm:Wor-Tan-Kre-thm}, the orbits of this action by $\Aut^\otimes(\Rep G)$ precisely correspond to the isomorphism classes of compact quantum groups with representation category $\Rep G$. 

\smallskip
Let us move on to $3$-cohomology. An invertible element $\Phi \in \U(G^3)$ is a $3$-cocycle if and only if it satisfies
$$
(1 \otimes \Phi)(\iota \otimes \hat{\Delta} \otimes \iota)(\Phi) (\Phi \otimes 1) = (\iota \otimes \iota \otimes \hat{\Delta})(\Phi) (\hat{\Delta} \otimes \iota \otimes \iota)(\Phi).
$$
Invariant $3$-cocycles are also called \emph{associators}. If $\Phi$ is such a cocycle, its action on $H_U \otimes H_V \otimes H_W$ can be considered as a new associativity morphism on the  C$^*$-category $\Rep G$ with bifunctor~$\circt$. If $\Phi$ is unitary, this gives a new C$^*$-tensor category $(\Rep G, \Phi)$, which has the same data as $\Rep G$ except for the new associativity morphisms defined by the action of~$\Phi$. If $\cC$ is a semisimple C$^*$-tensor category with the same fusion rules as $\Rep G$, that is, with the same Grothendieck semiring, then by transporting the monoidal structure of $\cC$ along any choice of a C$^*$-functor $\cC\to\Rep G$ defining the isomorphism $K^+(\cC)\cong R^+(G)$, we see that $\cC$ is monoidally equivalent to some $(\Rep G, \Phi)$ as above. We remark that in general it is not clear whether any such category is automatically rigid.

If $x$ is an invariant unitary $2$-cochain, the categories $(\Rep G, \Phi)$ and $(\Rep G, \Phi_x)$ are naturally unitarily monoidally equivalent, by means of the unitary tensor functor $$(\id_\un, \Id_{\Rep G}, x^{-1})\colon (\Rep G, \Phi) \to (\Rep G, \Phi_x).$$ This way the set $H^3_G(\hat{G}; \T)$ gives a parametrization of the categories of the form $(\Rep G, \Phi)$ considered up to unitary monoidal equivalences that preserve the isomorphism classes of objects.

\smallskip

A simple way of constructing elements of $H^3_G(\hat{G}; \T)$ is by considering the chain group of $\Rep G$, which we denote by $\Ch(G)$. Namely, any cocycle $\phi\in Z^3(\Ch(\cC);\T)$ can be considered as a $3$-cocycle on $\hat G$ such that its component in $B(H_r)\otimes B(H_s)\otimes B(H_t)$ is the scalar $\phi(g_r,g_s,g_t)$. We denote the category $\Rep G$ with the associator given by this cocycle by $(\Rep G)^\phi$. Therefore in $(\Rep G)^\phi$ the associativity morphism
$(U\circt V)\circt W\to U\circt (V\circt W)$ is the scalar operator $\phi(g,h,k)$ if $U$, $V$ and $W$ have the degrees $g$, $h$ and $k\in\Ch(G)$, respectively. The category $(\Rep G)^\phi$ is always rigid.

Let us note in passing that the Pontryagin dual of the abelianization $\Ch(G)^{\mathrm{ab}}$ of the chain group is naturally isomorphic to $H^1_G(\hat G;\T)$, cf.~\cite{MR2130607,arXiv:1506.09194}. Namely, any character $$\chi\colon\Ch(G)\to\T$$ defines a $1$-cocycle $(\chi(g_s))_{s\in\Irr\cC}\in\U(G)$, or equivalently, a monoidal automorphism $\eta^\chi$ of the identity functor on $\Rep G$ such that $\eta^\chi_U=\chi(g)\iota_U$ if~$U$ has degree $g\in\Ch(G)$.

In particular, if $G$ is a genuine group, then $\Ch(G)$ is abelian and we get an isomorphism $\Ch(G)\cong\widehat{Z(G)}$. More concretely, if $U$ is an irreducible representation of $G$, then the group~$Z(G)$ must be acting by scalars, that is, by some character $\chi_U \in \widehat{Z(G)}$, and the isomorphism $\Ch(G)\cong\widehat{Z(G)}$ maps $g_s$ into $\chi_{U_s}$.





\subsection{\texorpdfstring{$q$}{q}-Deformation}\label{sec:qdeform}

Let $G$ be a simply connected compact Lie group with a maximal torus~$T$, $\mathfrak{g}$ be the corresponding semisimple complex Lie algebra. Fix a system of simple roots and denote the Cartan matrix by $(a_{i j})_{i,j}$. We fix a scalar product defining the root system of $G$. When $G$ is simple, such a scalar product is unique up to a scalar factor and we normalize it so that for every short root $\alpha_i$ we have $(\alpha_i,\alpha_i)=2$.

For a fixed positive real number $q \neq 1$, the quantized universal enveloping algebra~$\U_q(\mathfrak{g})$ is the $*$-Hopf algebra (over $\C$) generated by elements $E_i, F_i, K_i^{\pm 1}$ satisfying
\begin{gather*}
[K_i, K_j] = 0, \quad K_i E_j K_i^{-1} = q_i^{a_{i j}} E_j, \quad K_i F_j K_i^{-1} = q_i^{-a_{i j}} F_j, \quad
[E_i, F_j] = \delta_{i j} \frac{K_i - K_i^{-1}}{q_i - q_i^{-1}},\\
\sum_{k = 0}^{1-a_{i j}} (-1)^k
\left [ \begin{array}{c}
  1 - a_{i j}\\
  k
\end{array} \right ]_{q_i}
E_i^k E_j E_i^{1 - a_{i j} - k} = 0,\\
\sum_{k = 0}^{1-a_{i j}} (-1)^k
\left [ \begin{array}{c}
  1 - a_{i j}\\
  k
\end{array} \right ]_{q_i}
F_i^k F_j F_i^{1 - a_{i j} - k} = 0.
\end{gather*}
where $q_i = q^{d_i}$ and $d_i=(\alpha_i,\alpha_i)/2$. The coproduct is
\begin{equation*}
\hat{\Delta}_q(E_i) = E_i \otimes 1 + K_i \otimes E_i,\quad
\hat{\Delta}_q(F_i) = F_i \otimes K_i^{-1} + 1 \otimes F_i, \quad
\hat{\Delta}_q(K_i) = K_i \otimes K_i,
\end{equation*}
and the involution is characterized by
$$
E_i^* = F_i K_i,\quad
F_i^* = K_i^{-1} E_i, \quad
K_i^* = K_i.
$$

A finite dimensional representation of $\U_q(\mathfrak{g})$ is said to be \emph{admissible} if it defines a representation of the maximal torus $T$ of $G$, so that the elements $K_i$ act as the elements $q_i^{H_i}$ in the complexification of $T$, where $H_i$ is the element of the Cartan subalgebra of $\mathfrak g$ defined by $\alpha_j(H_i)=a_{ij}$. The finite dimensional admissible unitary representations form a rigid C$^*$-tensor category. Since it is given as a subcategory of $\Hilbf$, Woronowicz's Tannaka--Krein theorem gives a compact quantum group, which is called the \emph{$q$-deformation} of $G_q$. The application of Woronowicz's theorem is a bit of an overkill here, and concretely, the Hopf $*$-algebra $\C[G_q]$ is defined as the subalgebra of~$\U_q(\mathfrak{g})^*$ generated by the matrix coefficients of admissible representations. It is known that the quantum group $G_q$ is of $G$-type according to Definition~\ref{def:type}.

The algebra of the discrete dual $\U(G_q)$ can be regarded as a completion of~$\U_q(\mathfrak{g})$ (although one should note that the latter has non-admissible representations which do not extend to $\U(G_q)$). The quantum group $G_q$ contains $T$ as a closed subgroup, with the corresponding embedding $\U(T)\hookrightarrow \U(G_q)$ given by identifying~$q_i^{H_i}$ with~$K_i$. Intrinsically the torus $T$ can be characterized as the \emph{maximal quantum Kac subgroup} of $G_q$, that is, the largest quantum subgroup of $G_q$ on which the antipode is involutive~\cite{MR2335776}.

The center $Z(G)\subset T$ of $G$ is contained in the center of $\U(G_q)$. It follows that for any subgroup $\Gamma\subset Z(G)$ the quantum group $G_q/\Gamma$ makes sense. This allows one to deform compact connected semisimple Lie groups that are not necessarily simply connected by letting $(G/\Gamma)_q=G_q/\Gamma$. In other words, $(G/\Gamma)_q$ is obtained from~$\U_q({\mathfrak g})$ by considering only the admissible representations such that their restrictions to $T$ factor through $T/\Gamma$, i.e., that have weights annihilating $\Gamma$.

\smallskip

It is known that for $q\ne1$ we have $H^1(\hat{G}_q;\T)=T$ and $H^1_{G_q}(\hat{G}_q;\T)=Z(G)$, which is a simple consequence of
So{\u\i}belman's classification of irreducible representations of $\C[G_q]$~\cite{MR1049910}. The second cohomology will be described in Sections~\ref{sec:cat-poisson} and~\ref{sec:functor-classification}. As for the third cohomology, since~$G_q$ and~$G$ are known to have the same fusion rules, there is a canonical bijection $H^3_{G_q}(\hat{G}_q; \T)\cong H^3_{G}(\hat{G}; \T)$. For the same reason there is a unique element in $H^3_{G}(\hat{G}; \T)$ corresponding to $\Rep G_q$. An associator~$\Phi_{\mathrm{KZ},q}$ representing this element was constructed by Drinfeld~\cite{MR1047964} using the Knizhnik--Zamolodchikov equations (see also~\cite{MR1239507,MR2832264}). Since the categories $\Rep G_q$ are mutually inequivalent for different $q\in(0,1]$, we thus get a family of different classes $[\Phi_{\mathrm{KZ},q}]\in H^3_{G}(\hat{G}; \T)$ indexed by $q\in(0,1]$. By our discussion in Section~\ref{sec:cohom-discr-dual}, there also are classes in $H^3_G(\hat{G}; \T)$ coming from $3$-cocycles on $\Ch(G)\cong\widehat{Z(G)}$, as well as from their products with~$\Phi_{\mathrm{KZ},q}$. What else is contained in $H^3_G(\hat{G}; \T)$, is a major open problem.


\section{Categorical duality for actions of quantum groups}
\label{sec:cat-duality-action}

In this section we describe several extensions of the Tannaka--Krein duality to actions of compact quantum groups. These extensions are not strictly speaking needed for the classification problem for compact quantum groups, but they provide motivation for some of the subsequent constructions.

\subsection{\texorpdfstring{$G$}{G}-algebras and \texorpdfstring{$(\Rep G)$}{Rep G}-module categories}

Given a compact quantum group $G$, a unital $G$-C$^*$-algebra is a unital C$^*$-algebra~$B$ equipped with a continuous left action $\alpha\colon B\to C(G)\otimes B$ of~$G$. This means that~$\alpha$ is an injective unital $*$-homomorphism such that $(\Delta\otimes\iota)\alpha=(\iota\otimes\alpha)\alpha$ and such that the space $(C(G)\otimes1)\alpha(B)$ is dense in $C(G)\otimes B$.
The linear span of spectral subspaces,
$$
\B = \{ x \in B \mid \alpha(x) \in \C[G] \otimes_\alg B \},
$$
which is a dense $*$-subalgebra of $B$, is called the \textit{regular subalgebra} of $B$, and the elements of~$\B$ are called \emph{regular}.  More concretely, the algebra $\B$ is spanned by the elements of the form $(h\otimes\iota)((x\otimes1)\alpha(a))$ for $x\in\C[G]$ and $a\in B$. This algebra is of central importance for the categorical reconstruction of~$B$.

Next let us explain the categorical counterpart of the $G$-algebras. We also note that a particular type of this structure plays a central role in the subfactor theory.

\begin{definition}
Let $\cC$ be a C$^*$-tensor category. We say that $\ccD$ is \emph{a right $\cC$-module category} if
\begin{itemize}
\item $\ccD$ is a C$^*$-category,
\item we are given a unitary bifunctor $\ccD \times \cC \to \ccD$, denoted by $(X, U) \mapsto X \times U$, and
\item natural unitary isomorphisms
\begin{equation}\label{eq:C-mod-untior-associator}
X \times \un \to X, \quad (X \times U) \times V \to X \times (U \otimes V)
\end{equation}
satisfying the obvious compatibility conditions analogous to those for monoidal categories.
\end{itemize}
\end{definition}

When $\ccD$ and $\ccD'$ are $\cC$-module categories, a \emph{$\cC$-module functor} $\ccD \to \ccD'$ is a pair $(\fF, \theta)$, where $\fF$ is a C$^*$-functor from $\ccD$ to $\ccD'$, and $\theta$ is a unitary natural transformation $\fF(X) \times U \to \fF(X \times U)$, again satisfying a standard set of compatibility conditions. When $\theta$ is obvious from context, we simply write $\fF$ instead of $(\fF, \theta)$. The $\cC$-module functors can be composed in a way similar to the monoidal functors.

We say that two $\cC$-module functors $(\fF, \theta), (\fF', \theta')\colon \ccD \to \ccD'$ are equivalent if there is a natural unitary transformation $\eta\colon \fF \to \fF'$ which is compatible with $\theta$ and $\theta'$. Note that this equivalence relation is compatible with composition of module functors.

Our starting point is the following categorical duality theorem for $G$-algebras, which could be called the Tannaka--Krein duality theorem for quantum group actions. Results leading to this theorem have a long history, starting from the work of Wassermann~\cite{MR0990110} and Landstad~\cite{MR1190512} in the early 1980s on full multiplicity ergodic actions of compact groups, and continuing in~\cite{MR1006625,MR1976459,MR2202309,MR2358289}.

\begin{theorem}[\cite{MR3121622,DOI:10.17879/58269764511}]
\label{tactions}
Let $G$ be a reduced compact quantum group.  Then the following two categories are equivalent:
\begin{enumerate}
\item The category of unital $G$-C$^*$-algebras, with unital $G$-equivariant $*$-homo\-mor\-phisms as morphisms.
\item The category of pairs $(\ccD,M)$, where $\ccD$ is a right $(\Rep G)$-module C$^*$-category and $M$ is a generating object in $\ccD$, with equivalence classes of unitary $(\Rep G)$-module functors respecting the prescribed generating objects as morphisms.
\end{enumerate}
\end{theorem}

The condition that $\ccD$ is generated by $M$ means that
any object in $\ccD$ is isomorphic to a subobject of $M\times U$ for some
$U\in\Rep G$.

\begin{remark}
Let $\End(\ccD)$ be the C$^*$-tensor category of C$^*$-endofunctors on $\ccD$, with uniformly bounded natural transformations as morphisms. Then, having a $\cC$-module structure on $\ccD$ is the same as giving a unitary tensor functor $\fF\colon \cC^{\otimes\mathrm{op}} \to \End(\ccD)$, where $\fF_0$ and $\fF_2$ correspond to the morphisms in~\eqref{eq:C-mod-untior-associator}. For $G = \SU_q(2)$ and, more generally, for free orthogonal quantum groups this point of view leads to a combinatorial classification of ergodic actions~\cite{MR3420332}.
\end{remark}

To describe the above equivalence, given a $G$-C$^*$-algebra $(B, \alpha)$, we consider the category~$\ccD_B$ of $G$-equivariant finitely generated right Hilbert $B$-modules.  In other words, objects of $\ccD_B$ are finitely generated right Hilbert $B$-modules $X$ equipped with a linear map $\delta=\delta_X\colon X\to C(G)\otimes X$ which satisfies the comultiplicativity property  $(\Delta\otimes\iota)\delta=(\iota\otimes\delta)\delta$ together with the following conditions:
\begin{itemize}
\item $(C(G)\otimes 1)\delta(X)$ is dense in $C(G)\otimes X$,
\item $\delta$ is compatible with the Hilbert $B$-module structure, in the sense that
\begin{equation*}
  \delta(\xi a) = \delta(\xi)\alpha(a), \quad
  \langle\delta(\xi),\delta(\zeta)\rangle = \alpha(\langle\xi,\zeta\rangle), \quad
  (\xi, \zeta\in X, a\in B).
\end{equation*}
\end{itemize}
Here, $C(G)\otimes X$ is considered as a right Hilbert $(C(G)\otimes B)$-module.

For $X \in \ccD_B$ and $U \in \Rep G$, we obtain a new object $X \times U$ in $\ccD_B$ given by the linear space $H_U \otimes X$, which is a right Hilbert $B$-module such that
$$
(\xi \otimes x) a = \xi \otimes x a, \quad \langle \xi \otimes x, \eta \otimes y \rangle_B = (\eta, \xi) \langle x, y \rangle_B \quad \text{for}\ \ \xi, \eta \in H_U,\ x, y \in X,\ a \in B,
$$
together with the compatible $C(G)$-coaction map
\begin{equation*}
\delta=\delta_{H_U\otimes X} \colon H_U\otimes X \to C(G)\otimes H_U\otimes X, \quad \delta(\xi\otimes x)=U^*_{21}(\xi\otimes\delta_X(x))_{213}.
\end{equation*}
This construction gives the structure of a right $(\Rep G)$-module category on $\ccD_B$, together with a distinguished object $B \in \ccD_B$.

\smallskip
In the other direction, we construct what would be the regular subalgebra of a $G$-algebra starting from a pair~$(\ccD,M)$ as in Theorem~\ref{tactions}. The generating condition on $M$ implies that, by replacing $\ccD$ by an equivalent category, we may assume that~$\ccD$ is the idempotent completion of the category $\Rep G$ with larger morphism sets $\ccD(U,V)$ than in $\Rep G$, such that~$M$ is the unit object $\un$ in $\Rep G$ and the functor $\iota\times U$ on $\ccD$ is an extension of the functor $\iota\circt U$ on $\Rep G$. Namely, we simply define the new set of morphisms between $U$ and $V$ as $\ccD(M\times U,M\times V)$.

Consider the linear space
\begin{equation*}
  \B=\bigoplus_{s\in\Irr G}\bar H_s\otimes \ccD(\un,U_s).
\end{equation*}
Note that this gives $\C[G]$ if $\ccD$ is $\Hilbf$ with the action of $\Rep G$ induced by the canonical fiber functor. Using this observation, the associative product on $\B$ is defined by
$$
(\bar\xi\otimes T)\cdot(\bar\zeta\otimes S)= \sum_{r\in\Irr G}  \overline{u^{r\alpha *}_{t s}(\xi\otimes\zeta)} \otimes u^{r\alpha *}_{t s} (T\otimes S),
$$
where $(u^{r\alpha}_{t s})_\alpha$ is an orthonormal basis in the space of morphisms $U_r \to U_t \circt U_s$.
Similarly, the $*$-structure on $\B$ is given by the following analogue of~\eqref{eq:star-on-CG}:
\begin{equation*}
(\bar\xi\otimes T)^* = \rho^{-1/2}\xi\otimes( T^* \otimes\iota)\bar R_t \quad (\bar{\xi} \otimes T \in\bar{H}_t\otimes \ccD(\un,U_t)).
\end{equation*}

The $*$-algebra $\B$ has a natural left $\C[G]$-comodule structure, defined by the map $\alpha\colon\B\to\C[G]\otimes\B$ such that if $\{\xi_i\}_i$ is an orthonormal basis in~$H_t$ and $u_{ij}$ are the matrix coefficients of~$U_t$ in this basis, then
\begin{equation*}
\alpha(\bar\xi_i\otimes T)=\sum_ju_{ij}\otimes\bar\xi_j\otimes T.
\end{equation*}
It is shown then that the action $\alpha$ is algebraic in the sense of \cite[Definition~4.2]{MR3121622}, meaning that the fixed point algebra $A=\B^G\cong\ccD(\un)$ is a unital C$^*$-algebra and the conditional expectation $(h\otimes\iota)\alpha\colon\B\to A$ is positive and faithful. It follows that there is a unique completion of~$\B$ to a C$^*$-algebra $B$ such that $\alpha$ extends to an action of the reduced form of $G$ on $B$. This finishes the construction of an action from a module category.

\subsection{Yetter--Drinfeld algebras and tensor functors}

An important class of $G$-C$^*$-algebras is the quantum homogeneous spaces defined by quantum subgroups. Namely, when $f\colon \C[G] \to \C[H]$ is a surjective homomorphism of Hopf $*$-algebras, the subalgebra
$$
\C[G/H] = \{ x \in \C[G]\colon (\iota \otimes f) \Delta(x) = x \otimes 1 \}
$$
completes to a $C(G)$-comodule subalgebra of $C(G)$. As the notation suggests, this should be regarded as the space of functions on the quotient $G/H$. The categorical counterpart to this structure is the $(\Rep G)$-module category $\Rep H$, where $\Rep G$ acts through the unitary tensor functor $\Rep G \to \Rep H$ induced by $f$.

In general, we could ask when a $(\Rep G)$-module category structure comes from a tensor functor, or in other words, when the $G$-equivariant Hilbert $B$-modules admit a tensor product operation. The answer is that the $G$-algebra structure of $B$ should be extended to that of a \emph{braided-commutative Yetter--Drinfeld $G$-algebra}. Let us briefly recall the relevant definitions.

Assume we have a continuous left action of $\alpha\colon B\to C(G)\otimes B$  of  a compact quantum group~$G$ on a unital C$^*$-algebra $B$, as well as a continuous right action $\beta\colon B\to \M(B\otimes c_0(\hat G))$ of the dual discrete quantum group~$\hat G$. The action $\beta$ defines a left $\C[G]$-module algebra structure $\rhd\colon \C[G]\otimes B\to B$ on~$B$ by
$$
x\rhd a=(\iota\otimes x)\beta(x)\ \ \text{for} \ \ x\in\C[G]\ \ \text{and}\ \ a\in B.
$$
Here we view $c_0(\hat G)$ as a subalgebra of $\U(G)=\C[G]^*$. This structure is compatible with involution, in the sense that
\begin{equation*}
x\rhd a^*=(S(x)^*\rhd a)^*.
\end{equation*}

\begin{definition}
We say that $B$ is a {\em Yetter--Drinfeld $G$-C$^*$-algebra} if the following identity holds for all $x\in\C[G]$ and $a\in\B$:
\begin{equation*}
\alpha(x\rhd a) =x_{(1)} a_{(1)}S(x_{(3)})\otimes ( x_{(2)}\rhd a_{(2)}),
\end{equation*}
where we use Sweedler's sumless notation, so we write the effect of $\Delta$ and $\alpha$ as $x_{(1)}\otimes x_{(2)}$, etc. A Yetter--Drinfeld $G$-C$^*$-algebra $B$ is said to be {\em braided-commutative} if for all $a,b\in\B$ we have
\begin{equation}\label{eBC}
ab=b_{(2)}(S^{-1}(b_{(1)})\rhd a).
\end{equation}

Note that when $b$ is in the fixed point algebra $A=\B^G$, the right hand side of the above identity reduces to $b a$, and we see that $A$ is contained in the center of $\B$.
\end{definition}


\begin{remark}
Yetter--Drinfeld $G$-C$^*$-algebras can be regarded as $D(G)$-C$^*$-algebras for the Drinfeld double $D(G)$ of $G$, and they are studied in the more general setting of locally compact quantum groups by Nest and Voigt~\cite{MR2566309}.
\end{remark}

The categorical counterpart of a braided-commutative Yetter--Drinfeld $G$-C$^*$-al\-gebra is a pair $(\cC,\fE)$, where $\cC$ is a C$^*$-tensor category and $\fE\colon\Rep G\to\cC$ is a unitary tensor functor such that~$\cC$ is generated by the image of $\fE$. The condition that~$\cC$ is generated by the image of $\fE$ means that
any object in $\cC$ is isomorphic to a subobject of $\fE(U)$ for some
$U\in\Rep G$. We stress that we do not assume that the
unit in $\cC$ is simple. In fact, the
C$^*$-algebra $\cC(\un)$ is exactly the fixed point algebra~$B^G$
in the C$^*$-algebra $B$ corresponding to~$(\cC,\fE)$ in the next theorem.

Define the morphisms  $(\cC,\fE) \to (\cC',\fE')$ as the equivalence classes of pairs  $(\fF,\eta)$, where~$\fF$ is a unitary tensor functor $\fF\colon\cC\to\cC'$ and $\eta$ is a natural unitary monoidal isomorphism $\eta \colon \fF\fE \to \fE'$. We say that $(\fF, \eta)$ and $(\fF', \eta')$ are equivalent if there is a natural unitary monoidal transformation of the unitary tensor functors $\xi\colon \fF \to \fF'$ which is compatible with~$\eta$ and~$\eta'$ in the sense that $\eta_U=\eta'_U\xi_{\fE(U)}\colon \fF(\fE(U))\to\fE'(U)$ for all $U\in\cC$. Again, this relation is compatible with the composition of morphisms, and we denote by $\Tens(\Rep G)$ the category of pairs $(\cC,\fE)$ with morphisms given by the equivalence classes of this relation.

\begin{theorem}[\cite{MR3291643}] \label{tcatch}
Let $G$ be a reduced compact quantum group.  Then the following two categories are equivalent:
\begin{enumerate}
\item\label{tcatch-item-alg}
The category $\YD(G)$ of unital braided-commutative Yetter--Drinfeld $G$-C$^*$-algebras, with unital $G$- and $\hat G$-equivariant $*$-homomorphisms as morphisms.
\item\label{tcatch-item-cat}
The category $\Tens(\Rep G)$.
\end{enumerate}
Moreover, given a morphism $[(\fF,\eta)]\colon (\cC,\fE)\to(\cC',\fE')$, the corresponding homomorphism of Yetter--Drinfeld C$^*$-algebras is injective if and only if $\fF$ is faithful, and it is surjective if and only if $\fF$ is~full.
\end{theorem}

To prove this theorem we have to enhance the ingredients of Theorem~\ref{tactions}. Suppose that $B$ is a Yetter--Drinfeld $G$-C$^*$-algebra. Then any $G$-equivariant right Hilbert $B$-module $X$ automatically admits a left action of $\B$, given by
$$
a \xi = \xi_{(2)} (S^{-1}(\xi_{(1)}) \rhd a) \quad (a \in \B,\ \xi \in X\ \text{such that}\ \delta_X(\xi)\in\C[G]\otimes_\alg X),
$$
that is, we interpret~\eqref{eBC} as a formula for the left action on $X$. When $B$ is braided-commutative, \eqref{eBC} guarantees that $X$ is a $\B$-bimodule, and $\ccD_B$ becomes a C$^*$-tensor category with tensor product given by $X\otimes Y=Y\otimes_B X$. Moreover, then $U \mapsto B \times U$ becomes a unitary tensor functor $\Rep G\to\ccD_B$.

\smallskip

In the opposite direction, suppose $\fE\colon\Rep G\to \cC$ is a unitary tensor functor.
Similarly to the previous subsection, we may assume that it is simply an embedding and construct a $G$-algebra $\B=\oplus_s(\bar H_s\otimes\cC(\un,U_s))$. We then have to define a $\C[G]$-module structure $\rhd\colon\C[G]\otimes\B\to\B$. It is defined using the map
$$
\tilde\rhd\colon (\bar H_s\otimes H_s)\otimes(\bar H_t\otimes\cC(\un,U_t))\to\overline{(H_s\otimes H_t\otimes\bar H_s)}\otimes\cC(\un,U_s\circt U_t\circt\bar U_s),
$$
\begin{equation*}
\label{eq:from-tensor-to-dual-action}
(\bar\xi\otimes\zeta)\trhd(\bar\eta\otimes T)=\overline{(\xi\otimes\eta\otimes\overline{\rho^{-1/2}\zeta})}\otimes(\iota\otimes T\otimes\iota)\bar R_s,
\end{equation*}
by decomposing $U_s\circt U_t\circt\bar U_s$ into irreducibles.

\smallskip

Let us now consider the case when the target category $\cC$ is $\Hilbf$. Then it can be shown that Theorem~\ref{tcatch} establishes a bijection between the isomorphism classes of unitary fiber functors $\Rep G\to \Hilb_f$ and the isomorphism classes of unital Yetter--Drinfeld $G$-C$^*$-algebras $B$ such that
\begin{itemize}
\item the $G$-algebra $\B$ is a Hopf--Galois extension of $\C$ over $\C[G]$, meaning that $\B^G=\C1$ and the Galois map
\begin{equation}\label{eq:Galois}
\Gamma\colon \B \otimes \B \to \C[G] \otimes \B, \quad x \otimes y \mapsto x_{(1)} \otimes x_{(2)} y,
\end{equation}
is bijective,
\item the $\C[G]$-module structure $\rhd\colon\C[G]\otimes\B\to\B$ is completely determined by the action of $G$ and coincides with the \emph{Miyashita--Ulbrich action}, defined by $$x \rhd a = \Gamma^{-1}(x\otimes1)_1 a \Gamma^{-1}(x\otimes1)_2.$$
\end{itemize}

In the purely algebraic setting the correspondence between fiber functors and Hopf--Galois extensions of $\C$ over $\C[G]$ was established in~\cite{MR1006625}. In the operator algebraic setting this was done in~\cite{MR2202309}. Note that in the last paper instead of bijectivity of the Galois map an equivalent condition of \emph{full quantum multiplicity} is considered.

\subsection{Dual category}\label{sec:dual-cat}

Two most natural braided-commutative Yetter--Drinfeld algebras associated with any compact quantum group $G$ are the algebras of functions on $G$ and on its discrete dual~$\hat G$. Let us consider the latter in more detail.

In the Hopf--von Neumann algebraic framework, the discrete dual is defined by
$$
\ell^\infty(\hat G)=\ell^\infty\mhyph\bigoplus_{s\in\Irr G} B(H_s)\subset\U(G).
$$
We have a left adjoint action $\alpha$ of $G$ on $\ell^\infty(\hat G)$ defined by
\begin{equation}\label{eadjoint}
B(H_s)\ni T\mapsto (U_s)_{21}^*(1\otimes T)(U_s)_{21}.
\end{equation}
This action is continuous only in the von Neumann algebraic sense, so in order to stay within the class of $G$-C$^*$-algebras, instead of $\ell^\infty(\hat G)$ we consider the norm closure $B(\hat G)$ of the regular subalgebra $\ell^\infty_{\alg}(\hat G)\subset\ell^\infty(\hat G)$. The right action $\Dhat$ of $\hat G$ on $\ell^\infty(\hat G)$ makes this algebra into a unital braided-commutative Yetter--Drinfeld C$^*$-algebra. In other words, the left $\C[G]$-module structure on~$\ell^\infty_{\alg}(\hat G)$ is defined by
\begin{equation*}\label{ead}
x\rhd a=(\iota\otimes x)\Dhat(a).
\end{equation*}


We want to describe the corresponding C$^*$-tensor category $\hat\cC=\cC_{B(\hat G)}$ and the unitary tensor functor $\fF=\fF_{B(\hat G)}\colon\Rep G\to\hat\cC$. The category $\hat\cC$ is the idempotent completion of the category with the same objects as in $\Rep G$, but with morphism sets
$\hat\cC(U,V)\subset B(H_U,H_V)\otimes\ell^\infty_{\alg}(\hat G)$. In this picture $\fF$ is the obvious embedding functor. For the reasons that will become apparent in a moment, it is more convenient to consider $\hat\cC(U,V)$ as a subset of $\ell^\infty_\alg(\hat G)\otimes B(H_U,H_V)$.
Then $\hat\cC(U,V)$ is the set of elements $T\in \ell^\infty_\alg(\hat G)\otimes B(H_U,H_V)$ such that
$$
V^*_{31}(\alpha\otimes\iota)(T)U_{31}=1\otimes T.
$$
From the definition of the adjoint action~$\alpha$ we see that an element $T\in\ell^\infty_\alg(\hat G)\otimes B(H_U,H_V)$ lies in $\hat\cC(U,V)$ if and only if it defines a $G$-equivariant map $H_s\otimes H_U\to H_s\otimes H_V$ for all~$s$. It follows that $\hat\cC(U,V)$ can be identified with the space of bounded natural transformations between the functors $\iota\circt U$ and $\iota\circt V$ on $\Rep G$,
$$
\Natb(\iota\circt U,\iota\circt V) \cong \ell^\infty\mhyph\bigoplus_{s\in\Irr G} \Hom_G(H_s \otimes H_U, H_s \otimes H_W).
$$
The composition of morphisms is defined in the obvious way as the composition of natural transformations. The tensor product of morphisms is determined by the tensor products $\iota_W\otimes\xi$ and $\xi\otimes\iota_W$ for $\xi=(\xi_X)_X\in \Natb(\iota\circt U,\iota\circt V)$. These are defined by
\begin{equation} \label{eq:tens-dual}
(\xi \otimes \iota_W)_X = \xi_X \otimes \iota_W, \quad (\iota_W \otimes \xi)_X = \xi_{X \circt W}.
\end{equation}

\section{Poisson boundary}

An essential ingredient in the classification of fiber functors is the notion of noncommutative Poisson boundary and its categorical counterpart, and the universality of the latter with respect to amenable tensor functors.

\subsection{Noncommutative Poisson boundary}

For a finite dimensional unitary representation $U$ of $G$, consider the state $\phi_U$ on~$B(H_U)$ defined~by
\begin{equation*} \label{ephi}
\phi_U(T) = \frac{\Tr(T \pi_U(\rho)^{-1})}{\dim_q U} \quad \text{for}\ \ T \in B(H).
\end{equation*}
If $U$ is irreducible, it can be characterized as the unique state satisfying
$$
(\iota \otimes \phi_U)(U_{2 1}^* (1 \otimes T) U_{2 1}) = \phi_U(T).
$$
For our fixed representatives of irreducible representations $\{U_s\}_s$ of $G$, we write $\phi_s$ instead of~$\phi_{U_s}$.

When $\phi$ is a normal state on $\ell^\infty(\hat{G})$, we define a completely positive map $P_\phi$ on~$\ell^\infty(\hat{G})$ by
$$
P_{\phi}(a) = (\phi \otimes \iota) \Dhat(a).
$$
If $\mu$ is a probability measure on $\Irr G$, we then define a normal unital completely positive map~$P_\mu$ on $\ell^\infty(\hat G)$ by $P_\mu = \sum_s \mu(s) P_{\phi_s}$.  The space
$$
H^\infty(\hat{G}; \mu) = \{ x \in \ell^\infty(\hat{G}) \mid x = P_\mu(x) \}
$$
of $P_\mu$-harmonic elements is called the \emph{noncommutative Poisson boundary}~\cite{MR1916370} of $\hat{G}$ with respect to $\mu$.  This is an operator subspace of $\ell^\infty(\hat{G})$ closed under the left adjoint action~$\alpha$ of~$G$ defined by \eqref{eadjoint} and the right action $\Dhat$ of $\hat G$ on itself by translations.  It has a new product structure
\begin{equation}\label{eq:choi-effros-prod}
x \cdot y = \lim_{n \to \infty} P_\mu^n(x y),
\end{equation}
where the limit is taken in the strong$^*$ operator topology. With this product $H^\infty(\hat{G}; \mu)$ becomes a von Neumann algebra (with the original operator space structure)~\cite{MR0430809}, which can be regarded as a generalization of the classical Poisson boundary~\cite{MR704539}. Moreover, the actions of~$G$ and~$\hat{G}$ on~$\ell^\infty(\hat G)$ define continuous actions on $H^\infty(\hat G; \mu)$ in the von Neumann algebraic sense.


\subsection{Categorical Poisson boundary}\label{sec:cat-poisson}

Using the categorical description of the discrete duals of compact quantum groups given in Section~\ref{sec:dual-cat}, we have a straightforward translation of noncommutative Poisson boundaries to its categorical counterpart~\cite{arXiv:1405.6572}.

Let $\cC$ be a strict rigid C$^*$-tensor category with simple unit. Consider the category~$\hat\cC$ defined as in Section~\ref{sec:dual-cat} for $\cC=\Rep G$. We have ``partial trace'' maps
$$
\tr_U \otimes \iota \colon \cC(U \otimes X, U \otimes Y) \to \cC(X, Y), \quad T \mapsto \frac{1}{d(U)} (R_U^* \otimes \iota) (\iota \otimes T) (R_U \otimes \iota).
$$
Using them we can define an operator $P_U$ on $\hat{\cC}(V, W)$ by
$$
(P_U(\eta))_X = (\tr_U \otimes \iota)(\eta_{U \otimes X}) \in \cC(X \otimes V, X \otimes W).
$$

As before, fix representatives $U_s$ of isomorphism classes of simple objects and write $P_s$ for~$P_{U_s}$.
Let $\mu$ be a probability measure on $\Irr \cC$. Define an operator $P_\mu$ acting on $\Natb(\iota \otimes V, \iota \otimes W)$ by $P_\mu=\sum_s \mu(s)P_{s}$. It is called the \emph{Markov operator associated with $\mu$}. A bounded natural transformation $\eta\colon\iota\otimes V\to\iota\otimes W$ is called $P_\mu$-{\em harmonic} if $P_\mu(\eta)=\eta$. We denote the set of $P_\mu$-harmonic  natural transformations $\iota\otimes V\to\iota\otimes W$ by $\cP(\cC; \mu)(V, W)$, or just by $\cP(V, W)$.

As is the case for usual harmonic functions, the naive product of harmonic transformation is not guaranteed to be harmonic. However, the formula~\eqref{eq:choi-effros-prod} still makes sense, and $\cP(\cC; \mu)(V, W)$ can be considered as a morphism set of a new C$^*$-category. We denote the subobject completion of this category by $\cP(\cC; \mu)$, or simply by $\cP$.

The category $\cP$ has the natural structure of a C$^*$-tensor category defined similarly to $\hat\cC$. Namely, at the level of objects (before the subobject completion), the monoidal product $\otimes$ is the same as in $\cC$, while the tensor product of morphisms is given by
$$
\xi\otimes\eta=(\xi\otimes\iota)\cdot(\iota\otimes\eta),
$$
where $\xi\otimes\iota$ and $\iota\otimes\eta$ are defined in the same way~\eqref{eq:tens-dual} as in $\hat\cC$.

Any morphism $T\colon V\to W$ in $\cC$ defines a bounded natural transformation $(\iota_X\otimes T)_X$ from $\iota\otimes V$ into $\iota\otimes W$, which is obviously $P_\mu$-harmonic for every $\mu$. This embedding of $\cC$-morphisms into $\cP$-morphisms defines a strict unitary tensor functor $\Pi\colon\cC\to\cP$.

\begin{definition}[\cite{arXiv:1405.6572}]
The pair $(\cP,\Pi )$ is called the \emph{Poisson boundary} of $(\cC,\mu)$.
\end{definition}

Note that $\hat\cC(\un)=\ell^\infty(\Irr \cC)$. A probability measure $\mu$ on $\Irr \cC$ is called \emph{ergodic}, if the only $P_\mu$-harmonic bounded functions on $\Irr \cC$ are the constant functions. This means that $\cP(\cC; \mu)$ still has a simple unit. It can be shown that such a $\mu$ exists if and only if the intrinsic dimension function on $\cC$ is \emph{weakly amenable} in the terminology of~\cite{MR1644299}, meaning that there exists a state on $\ell^\infty(\Irr \cC)$ which is $P_s$-invariant for all $s$. In this case we say that $\cC$ is \emph{weakly amenable}. For weakly amenable categories, the Poisson boundary has the following universal property.

\begin{theorem}[\cite{arXiv:1405.6572}]\label{thm:poisson}
Let $\cC$ be a rigid C$^*$-tensor category with simple unit, and~$\mu$ be an ergodic probability measure on $\Irr \cC$. Then the equality
$$
d^{\cP(\cC;\mu)}(\Pi(X)) = \norm{\Gamma_X}
$$
holds for all $X \in \cC$. Moreover, if $\cC'$ is another C$^*$-tensor category with simple unit, and $\fF\colon \cC \to \cC'$ is a unitary tensor functor such that $d^{\cC'}(\fF(X)) = \norm{\Gamma_X}$ holds for all~$X$, then there exists a unitary tensor functor $\fE\colon \cP(\cC;\mu) \to \cC'$, unique up to a natural unitary monoidal isomorphism, such that $\fE \Pi$ is naturally unitarily monoidally isomorphic to $\fF$.
\end{theorem}

Note that the subcategory generated by the image of $\fF$ is rigid, although the standard solutions for $\fF(U)$ might not be in the image of $\fF$.

\smallskip

The proof of the above theorem consists of two parts, which rely on very different techniques. One part establishes a universal property of the Poisson boundary among the functors defining the smallest dimension function on $\cC$. It relies on a study of certain completely positive maps and their multiplicative domains, which can be thought of as analogues of the classical Poisson integral. The second part shows that the smallest dimension function must be $X\mapsto\|\Gamma_X\|$. The proof relies on subfactor theory and is inspired by the works of Pimsner and Popa~\cite{MR1111570,MR1278111}. It is clear that the second part is not needed once we have at least one functor $\fF\colon\cC\to\cC'$ such that $d^{\cC'}(\fF(X)) = \norm{\Gamma_X}$, as is the case in the applications we describe below.

\smallskip

Let us draw some implications of Theorem~\ref{thm:poisson} for $\cC=\Rep G$ for a compact quantum group~$G$. Recall that $G$ is called \emph{coamenable} if $\hat G$ is amenable. By a Kesten-type criterion this is equivalent to requiring $\|\Gamma_U\|=\dim U$ for all $U\in\Rep G$. Furthermore, coamenability of $G$ implies weak amenability of $\Rep G$, so the above theorem can be applied. In our current setting the theorem says that there exists a universal unitary tensor functor $\Pi\colon\Rep G\to\cP$ such that $d^\cP(\Pi(U))=\dim U$ for all $U$. By universality and Woronowicz's Tannaka--Krein duality this functor must correspond to a quantum subgroup of $G$. Since the classical and quantum dimension functions coincide only in the Kac case, it is not difficult to figure out what this quantum subgroup must be and obtain the following result.

\begin{theorem}[\cite{DOI:10.1093/imrn/rnv241}] \label{thm:Kacuniverse}
Let $G$ be a coamenable compact quantum group, $K < G$ its maximal quantum subgroup of Kac type. Then the forgetful functor $\Rep G\to\Rep K$ is a universal unitary tensor functor defining the classical dimension function on $\Rep G$.
\end{theorem}

Using the correspondence between the noncommutative and categorical boundaries, this implies that if $G$ is a coamenable compact quantum group and $\mu$ is an ergodic probability measure on $\Irr G$, then the Poisson boundary $H^\infty(\hat{G}; \mu)$ is $G$- and $\hat G$-equivariantly isomorphic to~$L^\infty(G/K)$. This result was originally proved by Tomatsu~\cite{MR2335776} (he states the result in a more restricted form, but in fact his proof works in the generality we formulated).

\smallskip

What is more important for our deformation problems, Theorem~\ref{thm:Kacuniverse} sometimes allows us to reduce the classification of dimension-preserving fiber functors to an easier task.

\begin{corollary}[\cite{DOI:10.1093/imrn/rnv241}]\label{cor:coamen-2-cohom-classification-layman}
With $G$ and $K$ as in Theorem~\ref{thm:Kacuniverse}, there is a bijective correspondence between the (natural unitary monoidal) isomorphism classes of dimension-preserving unitary fiber functors $\Rep G\to\Hilbf$ and those of $\Rep K$. Namely, the correspondence maps a functor $\Rep K\to\Hilbf$ into its composition with the forgetful functor $\Rep G\to\Rep K$.

In other words, the natural map $H^2(\hat{K}; \T) \to H^2(\hat{G}; \T)$ induced by the inclusion $\U(K) \hookrightarrow \U(G)$ is a bijection.
\end{corollary}

This corollary is of course void of any content when $G$ is already of Kac type, e.g., when $G$ is a genuine compact group. In the latter case it is, however, still possible to say something interesting about $H^2(\hat G;\T)$ by the theory developed by Wassermann~\cite{MR0990110,MR1014926} and Landstad~\cite{MR1190512}. We will say a few words about it later.

\smallskip

Now, let us return to abstract C$^*$-tensor categories and give another example of a computation of the categorical Poisson boundary.

Recall that in Section~\ref{sec:cohom-discr-dual} we introduced twistings $(\Rep G)^\phi$ of $\Rep G$ by cocycles $\phi\in Z^3(\Ch(G);\T)$. The same construction makes sense for arbitrary $\cC$: given a cocycle $\phi\in Z^3(\Ch(\cC);\T)$, we define new associativity morphisms $(U\otimes V)\otimes W\to U\otimes (V\otimes W)$ as those of $\cC$ multiplied by the factor $\phi(g,h,k)$ if $U$, $V$ and $W$ have the degrees $g$, $h$ and $k\in\Ch(\cC)$, respectively, and denote the C$^*$-tensor category we thus obtain by $\cC^\phi$.

The category $\hat\cC$ is graded over $\Ch(\cC)$, since if $X$ and $Y$ are objects in $\cC$ of different degrees, then there are no nonzero morphisms $X\to Y$ in $\hat\cC$. It follows that for any probability measure $\mu$ on $\Irr\cC$ the Poisson boundary $\cP(\cC;\mu)$ is still graded over $\Ch(\cC)$. Therefore $\Ch(\cC)$ is a quotient of $\Ch(\cP(\cC;\mu))$, so any cocycle $\phi\in Z^3(\Ch(\cC);\T)$ can be viewed as a cocycle on $\Ch(\cP(\cC;\mu))$ and we can consider the twisted category $\cP(\cC;\mu)^\phi$. The functor $\Pi\colon\cC\to\cP(\cC;\mu)$ defines a tensor functor $\cC^\phi\to\cP(\cC;\mu)^\phi$, which we denote by $\Pi^\phi$. A natural question is how $(\cP(\cC;\mu)^\phi,\Pi^\phi)$ is related to the Poisson boundary of $(\cC^\phi,\mu)$. We have the following result.

\begin{theorem}[\cite{arXiv:1506.09194}]
Let $\cC$ be a rigid C$^*$-tensor category with simple unit, $\phi\in Z^3(\Ch(\cC);\T)$, and~$\mu$ be an ergodic probability measure on $\Irr \cC$. Then
$\Pi^\phi\colon\cC^\phi\to\cP(\cC;\mu)^\phi$ is a universal unitary tensor functor such that $d^{\cP(\cC;\mu)^\phi}(\Pi^\phi(X)) = \norm{\Gamma_X}$ for all~$X\in\cC^\phi$, hence it is isomorphic to the Poisson boundary of $(\cC^\phi,\mu)$.
\end{theorem}

This theorem is proved in \cite{arXiv:1506.09194} for $\cC=\Rep G$, but the general case is essentially the same. See also \cite{BNYdraft} for a related result.

From this we obtain the following useful extension of Theorem~\ref{thm:Kacuniverse}.

\begin{theorem}[\cite{arXiv:1506.09194}]\label{thm:pois-bdry-twist}
Let $G$ be a coamenable compact quantum group with maximal closed quantum subgroup $K$ of Kac type. Then
\begin{itemize}
\item $\Ch(G)$ is a quotient of $\Ch(K)$;
\item for any $3$-cocycle $\phi \in Z^3(\Ch(G);\T)$, the forgetful functor $(\Rep G)^\phi\to(\Rep K)^\phi$ is a universal unitary tensor functor defining the classical dimension function on~$(\Rep G)^\phi$.
\end{itemize}
\end{theorem}

\section{Quantum groups of Lie type}
\label{sec:q-grp-Lie-type}

\subsection{Fiber functors on twisted representation categories}\label{sec:functor-classification}

Let $G$ be a compact connected semisimple Lie group and $T\subset G$ be a maximal torus. We denote the weight and root lattices of~$\mathfrak g$ by~$P$ and~$Q$, respectively, and denote by~$X^*(T)\subset P$ the weight lattice of $T$. We also denote by $\Psi=(X^*(T),\Pi,X_*(T),\Pi^\vee)$ a fixed based root datum of $(G,T)$.

For any $q>0$ we have canonical isomorphisms
$$
\Ch(G_q)\cong\Ch(G)\cong\widehat{Z(G)}\cong X^*(T)/Q.
$$
As described in Section~\ref{sec:cohom-discr-dual}, for every $\phi\in Z^3(\widehat{Z(G)};\T)$ we can therefore consider a new category $(\Rep G_q)^\phi$ with associativity morphisms defined by the action of~$\phi$.

\begin{example}
Let us give concrete examples of such cocycles $\phi$. Assume $Z(G)$ is isomorphic to a cyclic group $\Z/n\Z$ (which happens, for example, when $G$ is simple and $G\not\cong\mathrm{Spin}(4n)$). Then any $3$-cocycle on $\Ch(G_q)\cong \Z/n\Z$ is cohomologous to the cocycle
$$
\phi(a, b, c) = \omega^{\left(\flr{\frac{a+b}{n}} - \flr{\frac{a}{n}} - \flr{\frac{b}{n}}\right) c} \quad (a, b, c \in \Z / n \Z)
$$
for some $n$-th root of unity $\omega$.
\end{example}

The first natural question is whether the categories $(\Rep G_q)^\phi$ admit any fiber functors.

Suppose that $c\in\U(T\times T)$ is a $\T$-valued $2$-cochain on the dual group $\hat{T} = X^*(T)$ such that its coboundary $\partial c$ is invariant under $Q$ in each variable. Then $\partial c$ can be considered a $3$-cocycle~$\Phi^c$ on~$\widehat{Z(G)}$ and we obtain the twisted category $(\Rep G_q)^{\Phi^c}=(\Rep G_q, \Phi^c)$ as above.
Since $\Phi^c$ is the coboundary of $c$ over $\hat{T}$, we have a unitary fiber functor $\fF_c\colon(\Rep G_q,\Phi^c)\to\Hilbf$ which is identical to the canonical fiber functor on $\Rep G_q$, except that the tensor structure is given~by
$$
(\fF_c)_2\colon H_U \otimes H_V \to H_{U \circt V},\ \ \xi \otimes \eta \mapsto c^*(\xi \otimes \eta).
$$

This functor defines a new compact quantum group~$G_q^c$ such that $\Rep G_q^c$ is unitarily mono\-idally equivalent to $(\Rep G_q, \Phi^c)$. Explicitly, similarly to the case of twisting by $2$-cocycles, $\C[G_q^c]=\C[G]$ as coalgebras, while the new $*$-algebra structure is defined by duality from $(\U(G_q),c\Dhat_q(\cdot)c^*)$.
%
Because $c$ is defined on $\hat{T}$, the coproduct of any element $a \in \U(T)$ computed in $\U(G_q^c)$ is the same as $\hat{\Delta}_q(a) = \hat{\Delta}(a)$.  In particular, $T$ is still a closed subgroup of $G_q^c$.

\begin{example}\label{ex:c-tau}
A concrete example of $c$ can be given as follows~\cite{MR3340190}. Let $r$ be the rank of $G$, and let $\tau = (\tau_1, \ldots, \tau_r) \in Z(G)^r$. Take any function $c_\tau\colon X^*(T) \times X^*(T) \to \T$ satisfying
$$
c_\tau(\lambda, \mu + Q) = f(\lambda, \mu), \quad c_\tau(\lambda + \alpha_i,\mu) = \langle\tau_i,\mu\rangle c_\tau(\lambda,\mu).
$$
Then $\Phi^{c_\tau}=\partial(c_\tau)$ is $Q$-invariant in each variable and hence can be considered as a $3$-cocycle on~$\widehat{Z(G)}$. Moreover, any other choice of $c_\tau$ only differs by a function on $(X^*(T)/Q)^2$ and defines exactly the same quantum group~$G_q^{c_\tau}$.

Let $T_\tau$ be the subgroup of $Z(G)$ generated by the components of $\tau$.  Then the algebra $\C[G_q^{c_\tau}]$ can be presented as a subalgebra of $\C[G_q] \rtimes \hat{T}_\tau$, which can be used to understand irreducible representations and $K$-theory of $C(G_q^{c_\tau})$~\cite{MR3340190}.

The construction of $G^{c_\tau}_q$ is a particular example of twisting by \emph{almost adjoint invariant cocentral actions} developed in~\cite{arXiv:1506.09194}.
\end{example}

As was already mentioned in Section~\ref{sec:qdeform}, for $q\ne1$, the torus $T$ is a maximal quantum subgroup of $G_q$ of Kac type. By Theorem~\ref{thm:pois-bdry-twist} it follows that for any $\phi\in Z^3(\widehat{Z(G)};\T) $ the forgetful functor
$(\Rep G_q)^\phi\to(\Rep T)^\phi$
is a universal unitary tensor functor defining the classical dimension function on $(\Rep G_q)^\phi$. This can be used in two ways. On the one hand, applying this to $\phi=\Phi^c$ we conclude that
\begin{equation}\label{eq:second-cohom}
H^2(\hat G^c_q;\T)\cong H^2(\hat T,\T).
\end{equation}
In particular, the dimension-preserving fiber functors of $G^c_q$ produce only quantum groups of the form $G^{c'}_q$ with $c'c^{-1}\in Z^2(\hat T;\T)$.
On the other hand, since it is easy to understand when $(\Rep T)^\phi$ admits a fiber functor, we conclude that $(\Rep G_q)^\phi$ admits a dimension-preserving unitary tensor functor if and only if $\phi$ lifts to a coboundary on $\hat T=X^*(T)$. This leads to the following result.

\begin{theorem}[\cite{MR3340190,arXiv:1506.09194}]
\label{thm:3coc-2pseudo}
Let $G$ be a compact connected semisimple Lie group with a maximal torus $T$, and $H$ be a compact quantum group of $G$-type with representation category $(\Rep G_q)^\phi$ for some $q>0$, $q\ne1$, and $\phi\in Z^3(\widehat{Z(G)};\T)$. Then
\begin{itemize}
\item the cocycle $\phi$ is cohomologous to a cocycle of the form $\Phi^{c_\tau}$ constructed in Example~\ref{ex:c-tau};
\item the quantum group $H$ is isomorphic to $G_q^{\theta c_\tau}$ for a cocycle $\theta\in Z^2(\hat T;\T)$.
\end{itemize}
\end{theorem}

In particular, although the construction in Example~\ref{ex:c-tau} may seem somewhat ad hoc, the theorem shows that the only other way of constructing dimension-preserving unitary fiber functors on $(\Rep G_q)^\phi$ is by multiplying the cochains $c_\tau$ in that example by $2$-cocycles on $\hat T$.

\begin{remark}\mbox{\ }

\noindent
1. Since the classical and quantum dimension functions on $(\Rep G_q)^\phi$ coincide if and only if $q=1$, instead of requiring $q\ne1$ we could say that $H$ is not of Kac type.

\smallskip\noindent
2. Another way of formulating the above theorem is by saying that for $q\ne1$ the quantum groups~$G_q^{\theta c_\tau}$ exhaust all quantum groups of $G$-type corresponding to the classes $[\phi\Phi_{\mathrm{KZ},q}]\in H^3_G(\hat G;\T)$.

\smallskip\noindent
3. The cocycles $\phi$ cohomologous to the cocycles of the form $\Phi^{c_\tau}$ can be abstractly characterized as those that vanish on $\wedge^3(X^*(T)/Q)\subset H_3(X^*(T)/Q;\Z)$~\cite{MR3340190}. In particular, when $G$ is simple, then $\wedge^3(X^*(T)/Q)=0$ and all cocycles satisfy this property.
\end{remark}

Isomorphism~\eqref{eq:second-cohom} and Theorem~\ref{thm:3coc-2pseudo} fail miserably for $q=1$. First of all, the map $H^2(\hat T;\T)\to H^2(\hat G;\T)$ maps the orbits under the action of the Weyl group into single points. Therefore as soon as the rank of some simple factor of $G$ is at least~$2$, the map $H^2(\hat T;\T)\to H^2(\hat G;\T)$ is not injective. More importantly, when the rank is large enough, the map is very far from being surjective. In order to formulate the precise result, let us introduce some terminology.

Let $G$ be for a moment an arbitrary compact group. A cocycle $x\in Z^2(\hat G;\T)$ is called \emph{nondegenerate} if its cohomology class does not arise from a proper closed subgroup of $G$. Denote by $H^2(\hat G;\T)^\times$ the (possibly empty) subset of $H^2(\hat G;\T)$ consisting of classes represented by nondegenerate cocycles. The following result is essentially due to Wassermann~\cite{MR0990110,MR1014926}, although the formulation is rather taken from~\cite{MR2844801}. For finite groups and without the unitarity condition similar results were also obtained in~\cite{MR1264314,MR1617921,MR1747109}.

\begin{theorem}
For any compact group $G$ we have a decomposition
$$
H^2(\hat G;\T)\cong\bigsqcup_{[K]}H^2(\hat K;\T)^\times/N_G(K),
$$
where the union is taken over the conjugacy classes of closed subgroups $K$ of $G$ and $N_G(K)$ denotes the normalizer of $K$ in $G$, which acts through the adjoint action on $H^2(\hat K;\T)$.
\end{theorem}

It can further be shown that a cocycle $x\in Z^2(\hat K;\T)$ is nondegenerate if and only if the corresponding $K$-C$^*$-algebra $B$, as described in Section~\ref{sec:cat-duality-action}, is simple (as a C$^*$-algebra)~\cite{MR0990110,MR1190512}. When $K$ is finite, the structure of such $K$-C$^*$-algebras is easy to understand: $B$ must be a full matrix algebra $\End(H)$, and by bijectivity of the Galois map~\eqref{eq:Galois} its dimension has to be~$|K|$, so the action of $K$ must be given by an irreducible projective representation $K\to PU(H)$ of dimension $\dim H=|K|^{1/2}$.

Now, given a compact connected Lie group $G$ with a maximal torus $T$, if $G$ contains a non-abelian finite subgroup $K$ admitting such a projective representation, the corresponding class in $H^2(\hat G;\T)$ does not lie in the image of $H^2(\hat T;\T)$. Furthermore, as there can be a lot of such nonconjugate subgroups $K$, the above theorem suggests that we should not expect a simple parametrization of Kac quantum groups of $G$-type similar to Theorem~\ref{thm:3coc-2pseudo}, in particular, those that have representation category $\Rep G$.

\subsection{Isomorphism problem for twisted \texorpdfstring{$q$}{q}-deformations}
\label{sec:autoequiv}

In this last section we consider the problem of classifying the quantum groups~$G^c_q$ up to isomorphism. Since we have already classified the dimension-preserving unitary fiber functors on $\Rep G^c_q$ (for $q\ne1$), for this we have to understand unitary monoidal equivalences between these categories, and in particular, their autoequivalences.

We start with autoequivalences that preserve the isomorphism classes of objects. It is not difficult to see that the group of such autoequivalences of $(\Rep G_q)^\phi$ does not depend on $\phi$~\cite{DOI:10.1093/imrn/rnv241}.
Therefore these autoequivalences are described by the following result.

\begin{theorem}[\cite{MR2844801,MR2959039}] \label{pinvcoh}
For any $q>0$ and any compact connected semisimple Lie group $G$, we have a group isomorphism
$$
H^2(\widehat{Z(G)}; \T)\cong H^2_{G_q}(\hat{G}_q; \T)
$$
induced by the inclusion $\U(Z(G)) \hookrightarrow \U(G_q)$.
\end{theorem}

In view of \eqref{eq:second-cohom} the result is not surprising, at least for $q\ne1$. The proof, however, relies on different ideas and is more constructive than that of~\eqref{eq:second-cohom}. The main part of it shows that if $x\in Z^2_{G_q}(\hat G_q;\T)$ is an invariant cocycle with the property that for any highest weight representations $U_\lambda$ and $U_\eta$ it acts trivially via the representation $U_\lambda\circt U_\eta$ on the vectors of weights $\lambda+\eta$ and $\lambda+\eta-\alpha_i$ for all simple roots $\alpha_i$, then $x=1$.

\smallskip

Assume now that we have a unitary monoidal equivalence between $(\Rep G_q)^{\phi_1}$ and $(\Rep G_q)^{\phi_2}$. It defines an automorphism of the representation semiring $R^+(G)$. Such an automorphism must arise from an automorphism $\sigma$ of the based root datum~$\Psi$ of~$G$~\cite{MR733774}. The automorphism $\sigma$ can then be lifted to an automorphism of $G_q$, which has to map the class of $\phi_1$ in $H^3_{G_q}(\hat G_q;\T)$ into that of $\phi_2$. At least for simple Lie groups it is not difficult to understand when this happens.

\begin{proposition}[\cite{MR1237835,DOI:10.1093/imrn/rnv241}]\label{prop:simple-G-3-cohom-inj}
For any $q>0$ and any compact connected simple Lie group $G$, the canonical map
$H^3(\widehat{Z(G)}; \T) \to H^3_{G_q}(\hat{G}_q;\T)$ is injective , or equivalently, the map
$$
H^3(\widehat{Z(G)}; \T) \to H^3_{G}(\hat{G};\T),\ \ [\phi]\mapsto[\phi\Phi_{\mathrm{KZ},q}],
$$
is injective. Furthermore, unless $G\cong\mathrm{Spin}(4n)$, the group $\Aut(\Psi)$ acts trivially on $H^3(\widehat{Z(G)}; \T)$.
\end{proposition}

The result is proved in~\cite{MR1237835} for $G = \SU(n)$. The proof for other simple groups is similar. The idea is that if $U$ is an irreducible representation of $G_q$ and $f\colon\un\to U^{\circt k}$ is an isometric morphism which up to a phase factor can be described entirely in terms of the fusion rules, then the composition
\begin{equation*}\label{eq:compos-f-assoc-f-star}
\xymatrix@C=2.5em{
U \ar[r]^-{f \otimes \iota} & (U \circt U^{\circt( k-1)}) \circt U \ar[r]^{\phi}& U \circt (U^{\circt ( k-1)} \circt U) \ar[r]^-{\iota \otimes f^*} & U
}
\end{equation*}
in $(\Rep G_q)^\phi$ is a scalar which, on the one hand, depends only on the class of $\phi$ in $H^3_{G_q}(\hat{G}_q;\T)$ and, on the other hand, can be explicitly computed in terms of $\phi$. It should be remarked that the proposition, as well as several other results in this section, is formulated in~\cite{DOI:10.1093/imrn/rnv241} for simply connected groups. But the proofs for their quotients are only easier. For example, the center of $\mathrm{Spin}(4n)$ is $\Z/2\Z\oplus\Z/2\Z$ and in the above argument several different representations~$U$ are needed to recover the class of $\phi$. But the center of any proper quotient of  $\mathrm{Spin}(4n)$ is cyclic and it suffices to take only one $U$, namely, either the standard representation of $\mathrm{SO}_q(4n)$ or one of the spin representations.

\smallskip

Returning to classification of the quantum groups $G^c_q$, we see that, at least for simple $G$, unitary monoidal equivalences between their representation categories arise only from automorphisms of the root datum and from $2$-cocycles on the dual of the center. Together with the classification of dimension-preserving fiber functors on $\Rep G^c_q$ given by \eqref{eq:second-cohom}, this leads to the following result.

\begin{theorem}[\cite{DOI:10.1093/imrn/rnv241}]\label{thm:simple-G-isom-class}
Let $G$ be a compact connected simple Lie group with a maximal torus $T$, $q_1,q_2\in(0,1)$, and $c_1,c_2$ be $\T$-valued $2$-cochains  on $X^*(T)$ such that $\partial c_1,\partial c_2$ descend to $X^*(T)/Q$. Then the quantum groups~$G^{c_1}_{q_1}$ and~$G^{c_2}_{q_2}$ are isomorphic if and only if $q_1=q_2$ and there exist an element $\sigma \in \Aut(\Psi)$ and a $\T$-valued $2$-cochain $b$ on $X^*(T)/Q$ such that $c_1 \sigma(c_2)^{-1}b^{-1}$ is a coboundary on $X^*(T)$.
\end{theorem}

If we knew that the categories $(\Rep G_q)^\phi$ exhausted all rigid C$^*$-tensor categories with fusion rules of $G$, then the above theorem together with Theorem~\ref{thm:3coc-2pseudo} would give a classification of all non-Kac compact quantum groups of $G$-type, giving a partial solution of Problem~\ref{prob:classification}. For $G=\SU(n)$ this is known to be the case~\cite{MR1237835,MR3266525} (see also~\cite{MR2307417,MR2825504} for related slightly weaker results). Furthermore, in this case the algebras $\C[\SU^c_q(n)]$ can be explicitly described in terms of generators and relations~\cite{DOI:10.1093/imrn/rnv241}.

\begin{remark}
There are indications that a complete classification of non-Kac compact quantum groups of $G$-type should be possible beyond the case of $G=\SU(n)$. Specifically, in the framework of rigid semisimple $\C$-linear tensor categories, there are several recent classification results showing that such tensor categories with fusion rules of $G$ are exhausted by $\Rep G_q$, with $q \in \C^\times$ not a nontrivial root of unity. For example, the case of $G = \mathrm{PSp}(4)$ is known by~\cite{arXiv:1410.2876},
and the case of the smallest exceptional group $G = \mathrm{G}_2$ has been recently settled in~\cite{arXiv:1501.06869}. For these cases it only remains to show the uniqueness of the  compatible C$^*$-structures when~$q$ is a positive real number.
\end{remark}

\providecommand{\bysame}{\leavevmode\hbox to3em{\hrulefill}\thinspace}
\providecommand{\MR}{\relax\ifhmode\unskip\space\fi MR }
\providecommand{\MRhref}[2]{%
  \href{http://www.ams.org/mathscinet-getitem?mr=#1}{#2}
}
\providecommand{\href}[2]{#2}

\bigskip

\end{document}